\documentclass[a4j,12pt]{article} 
\usepackage{amsmath,amssymb}
\usepackage{mathdots}   
\newcommand{\bfd}{{\bf d}}

\allowdisplaybreaks[1]   

\usepackage{tikz}
\usetikzlibrary{calc,decorations.markings}

\newcommand{\Section}[1]{%
\renewcommand{\thesection}{\S\arabic{section}}
\section{#1}
\renewcommand{\thesection}{\arabic{section}}
\setcounter{equation}{0}}

\newcommand{\qed}{\hfill \hbox{\rule[-2pt]{4pt}{7pt}}}
\newcommand{\proof}{{\hspace*{0.4cm} {\it Proof}.\ \enskip}}

\newtheorem{lem}{Lemma}[section]
\newtheorem{thm}[lem]{Theorem}  
\newtheorem{prop}[lem]{Proposition}
\newtheorem{cor}[lem]{Corollary}
\newtheorem{rem}[lem]{Remark}

\newcommand{\G}{{\Bbb G}}

\newcommand{\Z}{{\Bbb Z}}
\newcommand{\R}{{\Bbb R}}
\newcommand{\C}{{\Bbb C}}

\newcommand{\N}{{\Bbb N}}

\newcommand{\X}{{\Bbb X}}

\newcommand{\calC}{{\cal C}}
\newcommand{\calD}{{\cal D}}

\newcommand{\calH}{{\cal H}}

\newcommand{\calO}{{\cal O}}
\newcommand{\calP}{{\cal P}}
\newcommand{\calR}{{\cal R}}
\newcommand{\calS}{{\cal S}}

\newcommand{\frX}{{\frak X}}

\newcommand{\pair}[2]{\left\langle {#1}, \, {#2}\right\rangle}

\newcommand{\set}[2]{\left\{\left.#1\vphantom{#2}\:\right\vert\:#2\right\}}
\newcommand{\wt}{\widetilde}

\newcommand{\frp}{{\frak p}}

\newcommand{\lam}{{\lambda}}
\newcommand{\Lam}{{\Lambda}}

\newcommand{\alp}{{\alpha}}  
\newcommand{\vphi}{{\varphi}} 
\newcommand{\eps}{{\epsilon}} 

\newcommand{\slit}{\vspace{5mm}\noindent}
\newcommand{\mslit}{\vspace{3mm}\noindent}
\newcommand{\mmslit}{\vspace{1mm}\noindent}
\newcommand{\real}{{\rm Re}}

\newcommand{{\any}}{{}^\forall}
\newcommand{{\is}}{{}^\exists}
\newcommand{{\st}}{\; {\rm s.t.}\; }

\newcommand{\ol}[1]{\overline{#1}}

\newcommand{\hec}{{\calH(G,K)}}
\newcommand{\CKX}{{\calC^\infty(K \backslash X)}}
\newcommand{\SKX}{{\calS(K \backslash X)}}

\newcommand{\Llra}{{\Longleftrightarrow}}

\newcommand{\abs}[1]{\left\vert{#1}\right\vert}  

\newcommand{\tensor}{\mathop{\textstyle{\bigotimes}}}

\newcommand{\bcup}{\mathop{\textstyle{\bigcup}}}

\newcommand{\twomatrix}[4]{\begin{pmatrix}
                           {#1} & {#2}\\
                           {#3} & {#4}
                          \end{pmatrix}}

\newcommand{\mapright}[1]{\displaystyle{
   \smash{\mathop{\hbox to 1cm{\rightarrowfill}}^{#1}}}}

\newcommand{\mapdownl}[1]{\Big\downarrow
   \llap{$\vcenter{\hbox{$\scriptstyle#1\,$}}$ }}

\newcommand{\gyaddots}{%
\setlength{\unitlength}{1mm}
\begin{picture}(5,3.5)(-2,-1.5)
\put(0,0){$\cdot$}
\put(2,1.5){$\cdot$}
\put(-2,-1.5){$\cdot$}
\end{picture}}
\newcommand{\gyakuddots}{\smash{\lower0.3ex\hbox{\gyaddots}}}

\newcommand{\mapdownlr}[2]{\Big\downarrow
   \llap{$\vcenter{\hbox{$\scriptstyle#1$}}$ }
   \rlap{$\vcenter{\hbox{$\scriptstyle#2$}}$ }}

\begin{document}
\title{Spherical functions and local densities on the space of \\$p$-adic quaternion hermitian matrices
\footnote
{
This research is partially supported by Grant-in-Aid for Scientific Research JP16K05081. \hspace{17mm}
Keywords: Spherical functions, Plancherel formula, quaternion hermitian forms, Orthogonal polynomials.
Mathematical Subject Classification 2010: 11F85, 11E95, 11F70, 33D45, 33D52.
}}
\author{Yumiko Hironaka\\
}

\date{}
\maketitle



\setcounter{section}{-1}
\Section{Introduction}
Let $\G$ be a reductive algebraic group and $\X$ a $\G$-homogeneous affine algebraic variety, where everything is assumed to be defined over a $\frp$-adic field $k$. We denote by $G$ and $X$ the sets of $k$-rational points of $\G$ and $\X$, respectively. Taking a maximal compact subgroup $K$ of $G$, we consider the Hecke algebra $\hec$, that is the commutative $\C$-algebra generated by the characteristic functions of $KgK, \; g \in G$. Then, a nonzero $K$-invariant function on $X$ is called {\it a spherical function on $X$} if it is an $\hec$-common eigenfunction.

Spherical functions on homogeneous spaces comprise an interesting topic to investigate and a basic tool to study harmonic analysis on $G$-space $X$. 
Spherical functions on the spaces of sesquilinear forms  are particularly interesting, since  they can be regarded as generating functions of local densities of representations of such forms. For the cases of alternating forms of size $2n$ and unramified hermitian forms of size $n$, the main terms of the explicit formulas of spherical functions are related to Hall-Littlewood symmetric polynomials of type $A_n$, which are well studied, hence it is possible to extract local densities of forms (cf. \cite{HS1}, \cite{JMSJ}). 
For the space of unitary hermitian forms of size $m$, the main terms of the explicit formulas are related to Hall-Littlewood polynomials of type $C_n$, where $m=2n$ or $m=2n+1$, according to the parity of $m$, and the unitary group acting on $X$ is of type $C_n$ or $BC_n$, respectively (cf. \cite{HK1}, \cite{HK2}, \cite{Heven}).

In the present paper, we introduce the space $X$ of quaternion hermitian forms of size $n$ on a $\frp$-adic field $k$ and study spherical functions on it, where we assume $k$ has odd residual characteristic. In \S 1, we introduce Cartan decomposition of $X$ due to Jacobowitz and define typical spherical functions $\omega(x; s)$ on $X$. In \S 2, we introduce local densities of representations within quaternion hermitian forms, and give an induction theorem of spherical functions using local densities (Theorem~\ref{th: ind}). By this theorem, we may regard spherical functions as generating functions of local densities, and we give the explicit value of the local density of itself (Theorem~\ref{th: vol}). Then we define a spherical Fourier transform $F_0$ on the Schwartz space $\SKX$, which is an injective $\hec$-module map (Proposition~\ref{prop: sph tr}). In \S 3, we consider the functional equations and location of possible poles and zeros of $\omega(x; s)$ (Theorem~\ref{th: feq}). Then we introduce the normalized Fourier transform $F$ by modifying $F_0$, which gives an inclusion of $\SKX$ into the symmetric Laurent polynomial ring $\calR = \C[q^{\pm z_1}, \ldots, q^{\pm z_n}]^{S_n}$, where $\calR$ is isomorphic to $\hec$ by Satake transform (Theorem~\ref{th: spTrans}). In \S 4, we give the explicit formulas of $\omega(x; s)$ by a general method introduced in \cite{JMSJ},\cite{French} (Theorem~\ref{th: explicit}). In this case, we obtain a different kind of symmetric polynomials as the main terms of explicit formulas from those of other sesquilinear forms(cf. Remark~\ref{rem: comparison}). In \S 5, we study $\SKX$ more precisely for small $n$.  In \S 5.1,  for size $n \leq 4$, we determine the $\hec$-module structure of $\SKX$ and show the dimension of spherical functions on $X$ associated to general $z$ is equal to $1$. In \S 5.2, we introduce the Plancherel measure for size $2$ proved by Yasushi Komori and give the inversion formula.
\vspace{1cm}
\Section{The space $X$ and spherical functions on it}   
\newcommand{\calo}{{\mathfrak o}}
\newcommand{\Trd}{{\rm T_{rd}}}
\newcommand{\Nrd}{{\rm N_{rd}}}
\newcommand{\gen}[1]{\langle{#1}\rangle}

Let $k$ be a $\frp$-adic field, and denote by $\calo$ the ring of integers, $\pi$ a fixed prime element, $\frp = \pi\calo$, and $q$ the cardinality of $\calo/\frp$. Throughout this paper we assume $k$ has odd residual characteristic.  
Set $D$ be a division quaternion algebra over $k$, $\calO$ the maximal order in $D$, and $\calP$ the maximal ideal in $\calO$.  
Then there is an unramified quadratic extension $k'$ of $k$ in $D$, for which $k' = k(\eps), \; \eps^2 \in \calo^\times$ and we may take a prime element $\varPi$ of $D$ such that $\varPi^2 = \pi, \; \varPi \eps = -\eps \varPi$. Then the set $\{1, \eps, \varPi, \varPi \eps\}$ forms a basis for $\calO/\calo$ with the involution $*$ on $D$ defined by
\begin{eqnarray} \label{alp-alp*}
\alp = a+b\eps + c\varPi + d\varPi \eps \longmapsto \alp^* = a-b\eps - c\varPi - d\varPi \eps, \quad (a, b, c, d \in k),
\end{eqnarray}
and $\alp\alp^* \in k$.
There is a $k$-algebra inclusion $\vphi : D \longrightarrow M_2(k')$ such that
\begin{eqnarray}
&&
\vphi(\alp) = \twomatrix{a+b\eps}{(c-d\eps)\pi}{c+d\eps}{a-b\eps} \; \mbox{determined by}\; 
\alp(1 \; \varPi) = (1 \; \varPi)\vphi(\alp), \\
&&
\det(\vphi(\alp)) = \alp\alp^* = \Nrd(\alp) \in k,\\
&&
{\rm trace}(\vphi(\alp)) = \alp+\alp^* = \Trd(\alp) \in k,
\end{eqnarray}
where $\alp$ is written as in \eqref{alp-alp*}, $\Nrd$ is the reduced norm and $\Trd$ is the reduced trace on $D$.
Based on $\vphi$, we have a $k$-algebra inclusion $\vphi_n: M_n(D) \longrightarrow M_{2n}(k')$, and the reduced norm and the reduced trace of an element of $A \in M_n(D)$ give by
\begin{eqnarray}
&&
\Nrd(A) = \det(\vphi_n(A)), \; \Trd(A) = {\rm trace}(\vphi_n(A)) \; (\in k).
\end{eqnarray}
In particular, we see
\begin{eqnarray}
&&
\Nrd(a) = \det(a)^2, \; \Trd(a) = 2{\rm trace}(a), \quad \mbox{for \; } a \in M_n(k).
\end{eqnarray}
Since $\Nrd$ and $\Trd$ do not depend on the choice of splitting fields of $D$, we may use another $k$-algebra inclusion $\vphi_n': M_n(D) \longrightarrow M_{2n}(k(\varPi))$ based on 
\begin{eqnarray} \label{another inclusion of D}
&&
\vphi'(\alp) = \twomatrix{a+c\varPi}{(b+d\varPi)\eps^2}{b-d\varPi}{a-c\varPi} \in M_2(k(\varPi)), \quad
\alp(1 \; \eps) = (1 \; \eps)\vphi'(\alp), 
\end{eqnarray} 
One may refer for the above facts to Reiner's book [Re, \S 9, 13, 14]. 

We extend the involution $*$ on $A = (a_{ij}) \in M_{mn}(D)$ by $A^* = (a_{ji}^*) \in M_{nm}(D)$. 
We define the space $X_n$ of quaternion hermitian forms and the action of $G_n = GL_n(D)$ as follows
\begin{eqnarray}
&&
X_n = \set{A \in G_n}{A^* = A}, \\
&&
g \cdot A = gAg^* = A[g^*], \quad \mbox{for\; } (g, A) \in G_n \times X_n.
\end{eqnarray}
Denote by $K_n$ the maximal order in $G_n$, i.e., $K_n = G_n(\calO)$. Then,  it is known ([Jac, Theorem~6.2]) that the set $K_n \backslash X_n$ of $K_n$-orbits in $X_n$ is bijectively correspond to $\Lam_n$, where
\begin{eqnarray}
&&
\wt{\Lam_n} = \set{\alp = (\alp_1, \ldots, \alp_n) \in \Z^n}{\alp_1 \geq \alp_2 \geq \cdots \geq \alp_n}, \nonumber \\
&& \label{Lam-n}
\Lam_n = \set{\alp \in \wt{\Lam_n}}{
\mbox{if $\alp_i$ is odd, then $\sharp\set{j}{\alp_j = \alp_i}$ is even} }.
\end{eqnarray} 
In fact, writing $\alp \in \Lam_n$ as 
\begin{eqnarray} \label{shape of alp}
\alp  = (\underbrace{\gamma_1,\ldots, \gamma_1}_{m_1},\ldots, \underbrace{\gamma_r,\ldots, \gamma_r}_{m_r}), \; \gamma_1>\cdots >\gamma_r, \; m_j > 0, \; \sum_j m_j = n, 
\end{eqnarray}
one may take the matrix $\pi^\alp = \gen{\pi^{\gamma_1^{m_1}}} \bot \cdots \bot \gen{\pi^{\gamma_r^{m_r}}} \in X_n$, where 
\begin{eqnarray}
\gen{\pi^{\gamma^m}} = \left\{\begin{array}{ll}
Diag(\pi^e, \ldots, \pi^e) & \mbox{if}\; \gamma = 2e\\
\twomatrix{0}{\pi^e\varPi}{-\pi^e\varPi}{0} \bot \cdots \bot \twomatrix{0}{\pi^e\varPi}{-\pi^e\varPi}{0} 
& \mbox{if}\; \gamma = 2e+1
\end{array} \right\} \in X_m.  \nonumber
\end{eqnarray}
Set $\Lam_n^+ = \set{\alp \in \Lam_n}{\alp_n \geq 0}$ and $X_n^+ = X_n \cap M_n(\calO)$. Then
\begin{eqnarray}   \label{integral X}
X_n^+ = \cup_{\alp \in \Lam_n^+}\, K_n \cdot \pi^\alp.
\end{eqnarray} 
It is easy to see 
\begin{eqnarray} \label{Nrd-alp}
&&
\Nrd(\pi^\alp) = \pi^{\abs{\alp}}, \quad \abs{\alp} = \sum_{i=1}^n \alp_i \in 2\Z, \quad (\alp \in \Lam_n)
\end{eqnarray}
hence we have
\begin{eqnarray}
\Nrd(x) \in k^2, \quad \mbox{for }x \in X_n.
\end{eqnarray} 
For $g \in G = G_n$, we denote by $g^{(i)}$ the upper left $i \times i$-block of $g$, \; $1 \leq i \leq n$.  
We take the Borel subgroup $B = B_n$ of $G$ consisting of lower triangular matrices. Then for $(p, x) \in B \times X_n$, we have 
\begin{eqnarray}
\Nrd((p\cdot x)^{(i)}) &=&\Nrd(p^{(i)}\cdot x^{(i)}) \nonumber\\
&=& \psi_i(p)^2 \Nrd(x^{(i)}), \quad  \psi_i(p) = \Nrd(p^{(i)}), \; 1 \leq i \leq n.
\end{eqnarray}
Thus, for $x \in X=X_n$, we may define $d_i(x) \in k$ by $d_i(x)^2 = \Nrd(x^{(i)}), \; 1 \leq i \leq n$. Then, $d_i(x)$ is a $B$-relative invariant associated with $k$-rational character $\psi_i, \; 1 \leq i \leq n$.  
For $x \in X$ and $s \in \C^n$, we consider the integral
\begin{eqnarray}   \label{def: sph}
\omega(x; s) = \int_{K_n} \abs{\bfd(k\cdot x)}^s dk, \quad \abs{\bfd(y)}^s = \left\{\begin{array}{ll}
\prod_{i=1}^n \abs{d_i(y)}^{s_i} & \mbox{if \; } y \in X^{op}\\
0 & \mbox{otherwise},
\end{array} \right.
\end{eqnarray}
where $dk$ is the normalized Haar measure on $K = K_n$, $\abs{\; }$ is the absolute value on $k$ and
\begin{eqnarray}
X_n^{op} = \set{x \in X_n}{d_i(x) \ne 0, \; \mbox{for all}\; 1 \leq i \leq n}.
\end{eqnarray}
The integral in \eqref{def: sph} is absolutely convergent if $\real(s_i) \geq 0, \; 1 \leq i \leq n-1$, and continued to a rational function of $q^{s_1}, \ldots, q^{s_n}$ (cf. \cite[Remark~1.1]{JMSJ}), where $s_n$ is free because $\abs{d_n(k\cdot x)} = \abs{x}$ for $k \in K_n$. 
Then it becomes an element of 
\begin{eqnarray}   \label{CKX}
\CKX  = \set{\Psi: X \longrightarrow \C}{\Psi(k\cdot x) = \Psi(x), \; k \in K},
\end{eqnarray}
and we use the notation $\omega(x; s)$ in such sense.  
Denote by $\hec$ the Hecke algebra of $G$ with respect to $K$. We recall the action of $\hec$ on $\CKX$: 
\begin{eqnarray}  \label{hec acts on CKX}
f* \Psi (x) = \int_G f(g)\Psi(g^{-1}\cdot x)dg, \quad (f \in \hec, \; \Psi \in \CKX, \; x \in X),
\end{eqnarray}
where $dg$ is the normalized Haar measure on $G$. 
We call $\omega(x; s)$ {\it a spherical function on $X$}, since  
it is a common eigenfunction with respect to the above action of $\hec$, in fact 
\begin{eqnarray}
(f * \omega(\; ;s))(x) 
& =&  \lam_s(f) \omega(x; s),  \quad (f \in \hec).
\end{eqnarray}
Here $\lam_s$ is the $\C$-algebra map 
\begin{eqnarray}
&\lam_s: & \hec \longrightarrow \C(q^{s_1},\ldots, q^{s_n}), \nonumber\\
&&  f \longmapsto \int_{B}f(p)\prod_{i=1}^n \abs{\psi_i(p)}^{-s_i}\delta(p)dp,   \label{lam_s}
\end{eqnarray}
where $dp$ is the left invariant measure on $B_n$ with modulus character $\delta$. The Weyl group $S_n$ of $G$ acts on $\{s_1, \ldots, s_n\}$ through its action on the rational characters $\set{\abs{\psi_i}^{s_i}}{1 \leq i \leq n}$. It is convenient to introduce a new variable $z \in \C^n$ related to $s \in \C^n$ by
\begin{eqnarray} \label{change of var}
s_i = -z_i+z_{i+1}-2 \; \; (1 \leq i \leq n-1), \quad s_n = -z_n+n-1,
\end{eqnarray} 
and denote $\omega(x; s) = \omega(x; z)$ and $\lam_s = \lam_z$.  Then $S_n$ acts on $\{z_1, \ldots, z_n\}$ by permutation of indices, and the $\C$-algebra map $\lam_z$ is the Satake isomorphism
\begin{eqnarray}  \label{lam_z}
\lam_z: \hec \stackrel{\sim}{\longrightarrow} \C[q^{\pm z_1}, \ldots, q^{\pm z_n}]^{S_n}.
\end{eqnarray}
Because of this isomorphism, all the spherical functions on $X$ are parametrized by eigenvalues $z \in \C^n$ through $\hec \longrightarrow \C, \; f \longmapsto \lam_z(f)$, and $\lam_z$ is determined by the class of $z$ in $\left(\C\big{/}\frac{2\pi\sqrt{-1}}{\log q}\Z\right)^n\big{/}S_n$.

\vspace{1cm}

\Section{Local densities and spherical functions}   

{\bf 2.1.} We will give the induction theorem (Theorem~\ref{th: ind}) of spherical functions by means of local densities, by which we may regard spherical functions as generating functions of local densities of representations.  We start with the definition of local densities.
For $A \in X_m^+$  and $B \in X_n^+$ with $m \geq n$, we define the local density $\mu(B, A)$ and primitive local density $\mu^{pr}(B, A)$ of $B$ by $A$ as follows:
\begin{eqnarray}
&& 
\mu(B, A) = \lim_{\ell \rightarrow \infty} \frac{N_\ell(B, A)}{q^{\ell n(4m-2n+1) + n(n-1)}},\nonumber\\
&& \label{def: local density}
\mu^{pr}(B, A) = \lim_{\ell \rightarrow \infty} \frac{N_\ell^{pr}(B, A)}{q^{\ell n(4m-2n+1) + n(n-1)}}.
\end{eqnarray}
Here 
\begin{eqnarray}
&&
N_\ell(B, A) = \sharp \set{\ol{u} \in M_{mn}(\calO)\big{/}M_{mn}(\calP^{2\ell})}{A[u] - B \in M_n(\calP^{2\ell-1})}, \nonumber \\
&& \label{def: N(B,A)}
N_\ell^{pr}(B, A) = \sharp\set{\ol{u}  \in M_{mn}^{pr}(\calO/\calP^{2\ell})}{A[u] - B \in M_n(\calP^{2\ell-1})},
\end{eqnarray}
where we identify $M_{mn}(\calO)\big{/}M_{mn}(\calP^{2\ell})$ with $M_{m n}(\calO/\calP^{2\ell})$ and denote by $\ol{u}$ its element represented by $u \in M_{mn}(\calO)$, and an element in $M_{mn}(\calO)$ is called {\it primitive} if it belongs to the set $GL_m(\calO)\begin{pmatrix}{1_n}\\{0}\end{pmatrix}$, and we write 
\begin{eqnarray}  \label{primitive set}
M_{mn}^{pr}(\calO/\calP^{2\ell}) = \set{\ol{u} \in M_{mn}(\calO/\calP^{2\ell})}{u \mbox{ is primitive }}.
\end{eqnarray}
%
%
\begin{rem}{\rm 
The above definition is well-defined, since the conditions in \eqref{def: N(B,A)} and \eqref{primitive set} are
independent of the choice of the representative $u$ of $\ol{u}$. 
If $\ell$ is sufficiently large, then the $K_n$-orbit $K_n \cdot B$ decomposes into a finite union of the set $B_i + M_n(\calP^{2\ell})$, and the ratios in the right hand sides of \eqref{def: local density} becomes stable. This phenomenon is characteristic of local densities for sesquilinear forms (cf. [Ki2], [HS1, \S 3], [H1, \S 2]
). 
}\end{rem} 

\medskip
We note that, for a matrix $C = C^* \in M_n(D)$, $C$ belongs to $M_n(\calP^{2\ell-1})$ if and only if $C$ belongs to $H_n(\calP, \ell)$, where
\begin{eqnarray}
H_n(\calP, \ell) = \set{A = (a_{ij}) \in M_n(\calO)}{A = A^*, \; a_{ii} \in \frp^\ell, \; a_{ij} \in \calP^{2\ell-1}, \; (\any i, j)}.
\end{eqnarray}
%
By definition, $\omega(x; s)$ takes the same value on the $K_n$-orbit containing $x \in X_n$, further we see that 
\begin{eqnarray*}
\omega(\pi^r x; s) = q^{-r\sum_{i=1}^n is_i} \omega(x;s) = q^{r(z_1+\cdots +z_n)}\omega(x;s), \quad r \in \Z.
\end{eqnarray*}
Hence it suffices to show the induction theorem for $\pi^\xi, \, \xi \in \Lam_m^+$. 

\begin{thm} \label{th: ind} 
Let $m > n$ and assume that $\real(s_i) \geq 0$ for any $1 \leq i \leq n$. Then, for any $\xi \in \Lam_m^+$, one has
\begin{eqnarray}
\lefteqn{\omega(\pi^\xi; s_1, \ldots, s_n, 0, \ldots, 0)} \nonumber\\
&=&
\frac{w_n(q^{-2}) w_{m-n}(q^{-2})}{w_m(q^{-2})} \times \sum_{\alp \in \Lam_n^+} \frac{\mu^{pr}(\pi^\alp, \pi^\xi)}{\mu(\pi^\alp, \pi^\alp)} \cdot \omega(\pi^\alp; s_1\ldots, s_n) \nonumber \\
&=&
\frac{w_n(q^{-2}) w_{m-n}(q^{-2})}{w_m(q^{-2})} \prod_{i=1}^n(1-q^{-(s_i+ \cdots +s_n+2m-2i+2)}) \times \sum_{\alp \in \Lam_n^+} \frac{\mu(\pi^\alp, \pi^\xi)}{\mu(\pi^\alp, \pi^\alp)} \cdot \omega(\pi^\alp; s_1\ldots, s_n), \nonumber 
\end{eqnarray}
where $w_m(t) = \prod_{i=1}^m(1-t^i)$.
\end{thm}

\medskip
The above theorem can be proved in a similar way to the case for the other sesquilinear forms, i.e. alternating, hermitian and symmetric forms, so we omit the proof (cf. [HS1, Theorem~5], [H1, \S2 Theorem]). For the present case the result is proved in the master thesis of Y.~Ohtaka ([OY]) in a slightly different definition, and he used it to study the explicit formula of spherical functions of size 2. 

\medskip
The density $\mu(\pi^\alp, \pi^\alp) = \mu^{pr}(\pi^\alp, \pi^\alp)$ is given as follows, which we will prove in \S 2.2.  In \S 2.3, we will introduce a spherical transform $F_0$ on the Schwartz space on $X$ and show it is injective by using Theorem~\ref{th: ind} (Proposition~\ref{prop: sph tr}).

\begin{thm} \label{th: vol}
Assume $\alp \in \Lam_n$ is given as in \eqref{shape of alp}.
Then one has
\begin{eqnarray} \label{density oneself}
\mu(\pi^\alp, \pi^\alp) = q^{2n(\alp) + \frac12\abs{\alp} + \frac12\sharp\set{i}{\alp_i \mbox{ is odd}} } 
\prod_{j=1}^r\, \left\{\begin{array}{ll} 
w_{m_j}(-q^{-1}) & \mbox{if } 2 \mid \gamma_j\\[2mm]
w_{\frac{m_j}{2}}(q^{-4}) & \mbox{if } 2 \not{\mid}\, \gamma_j
\end{array} \right\}, 
\end{eqnarray}
where
\begin{eqnarray*} \label{w(t)}
n(\alp) = \sum_{i=1}^n(i-1)\alp_i, \qquad \abs{\alp} = \sum_{i=1}^n \alp_i.
\end{eqnarray*}
\end{thm}

\slit
{\bf 2.2.} In the following, ${}^{(pr)}$ means that the identity holds with and without the condition primitive, respectively.

\begin{prop} \label{lem: shift}
For $A \in X_m^+$ and $B \in X_n^+$ with $m \geq n$ and $e \in \N$, one has
\begin{eqnarray}
\mu^{(pr)}(\pi^e B, \pi^e A) = q^{en(2n-1)}\mu^{(pr)}(B, A).
\end{eqnarray}
\end{prop}

\proof
Assume $\ell$ is sufficiently large, and take $X \in M_{m n}(\calO)$ such that $A[X]- B \in H_n(\calP, \ell)$. For any $Y \in M_{m n}(\calO)$, one has
\begin{eqnarray*}
(\pi^e A)[X+\pi^\ell Y] - \pi^eB &=& 
\pi^e(A[X]-B) + \pi^{e + \ell}(Y^*AX + X^*AY + \pi^\ell YAY^*) \\
&\in &H_n(\calP, e+\ell),
\end{eqnarray*}
and $X + \pi^\ell Y$ is primitive if $X$ is. Hence $N_{e+\ell}^{(pr)}(\pi^e B, \pi^e A) = q^{4emn} N_\ell^{(pr)}(B, A)$, and 
\begin{eqnarray*}
\mu^{(pr)}(\pi^e B, \pi^e A) &=& 
\lim_{\ell \rightarrow \infty} \frac{N_{e + \ell}^{(pr)}(\pi^e B, \pi^e A)}{q^{(e+\ell)n(4m-2n+1)+n(n-1)}} \\
&=&
\lim_{\ell \rightarrow \infty} \frac{N_\ell^{(pr)}(B, A) q^{4emn}}{q^{\ell n(4m-2n+1)+n(n-1) + 4emn-en(2n-1)}} \\
&=&
q^{en(2n-1)}\mu^{(pr)}(B, A).
\end{eqnarray*} 
\qed

\begin{rem} \label{rem: shift}
{\rm 
Owing to Proposition~\ref{lem: shift}, we may define local density and primitive local density for any $A \in X_m$ and $B \in X_n$ with $m \geq n$ as follows:
Taking $e \in \N$ for which $\pi^e A \in X_m^+$ and $\pi^e B \in X_n^+$, 
\begin{eqnarray}
\mu^{(pr)}(B, A) = q^{-en(2n-1)}\mu^{(pr)}(\pi^eB, \pi^eA).
\end{eqnarray}
Then, wee see that Proposition~\ref{lem: shift} is valid for any $A \in X_m$, $B \in X_n$, and $e \in \Z$. 
}\end{rem}

\begin{prop} \label{lem: decomp}
Assume $\alp \in \Lam_m^+$ is decomposed as $\alp = (\gamma, \beta)$ with $\beta \in \Lam_n^+$ and $\gamma \in \Lam_{m-n}^+$. 
Then
\begin{eqnarray}
\mu(\pi^\alp, \pi^\alp) =  q^{2(m-n)\abs{\beta}} \mu^{pr}(\pi^\beta, \pi^\alp) \mu(\pi^\gamma, \pi^\gamma).
\end{eqnarray}
In particular, if $\gamma_{m-n} > \beta_1$, then $\mu^{pr}(\pi^\beta, \pi^\alp) = \mu(\pi^\beta, \pi^\beta)$ and 
\begin{eqnarray} \label{eq2 in lem}
\mu(\pi^\alp, \pi^\alp) =  q^{2(m-n)\abs{\beta}} \mu(\pi^\beta, \pi^\beta) \mu(\pi^\gamma, \pi^\gamma).
\end{eqnarray}
\end{prop}

\proof
We use the notation $\check{\pi}^\alp = j_m \cdot \pi^\alp$, where $j_m$ is the matrix of size $m$ such that all the anti-diagonal entries are $1$ and other entries are $0$. Then $\check{\pi}^\alp = \twomatrix{\check{\pi}^\beta}{0}{0}{\check{\pi}^\gamma}$, where $\check{\pi}^\beta$ and $\check{\pi}^\gamma$ are defined similarly.
Assume $\ell$ is sufficiently large, and take 
\begin{eqnarray} \label{fixed X}
\ol{X} \in M_{mn}^{pr}(\calO/\calP^{2\ell}) \; \mbox{such that} \; \pi^\alp[X] - \check{\pi}^\beta \in {H_n(\calP, \ell)}. 
\end{eqnarray}
For an extension $Y = (X Z) \in GL_m(\calO)$ of $X$, we have 
\begin{eqnarray} \label{cont Y}
\pi^\alp[Y] = \twomatrix{\pi^\alp[X]}{X^*\pi^\alp Z}{Z^*\pi^\alp X}{\pi^\alp[Z]}, \quad  \pi^\alp[X] - \check{\pi}^\beta \in {H_n(\calP, \ell)},
\end{eqnarray}
and we may assume that $X^*\pi^\alp Z \equiv 0 \pmod{\calP^{2\ell}}$ after changing the extension (since $\beta_1 \leq \gamma_{m-n}$), then $\pi^\alp[Z]$ is $K_{m-n}$-equivalent to $\check{\pi}^\gamma$.  Hence there is an extension $Y$ of $X$ such that 
\begin{eqnarray} \label{second cont Y}
\pi^\alp[Y] - \check{\pi}^\alp \in {H_m(\calP,\ell)},
\end{eqnarray}
or equivalently 
\begin{eqnarray} \label{third cont Y}
\check{\pi}^\alp[Y^{-1}] - \pi^\alp  \in {H_m(\calP,\ell)}.
\end{eqnarray}
For such extensions $Y_1$ and $Y_2$ of $X$, we see 
\begin{eqnarray}
\check{\pi}^\alp[Y_1^{-1}Y_2] - \check{\pi}^\alp \in H_m(\calP, \ell), \quad
Y_2 = Y_1 \twomatrix{1_n}{W}{0}{V} (\mbox{ in } GL_m(\calO)).
\end{eqnarray}
Since
\begin{eqnarray*}
\check{\pi}^\alp\left[\twomatrix{1_n}{W}{0}{V}\right] = \twomatrix{\check{\pi}^\beta}{\check{\pi}^\beta W} {W^*\check{\pi}^\beta}{\check{\pi}^\beta [W] + \check{\pi}^\gamma[V]}, 
\end{eqnarray*}
and $\ell$ is large enough, we see the number of extensions $Y$ of type \eqref{second cont Y} for the fixed $X$ as in \eqref{fixed X} is equal to 
\begin{eqnarray*}
&&
\sharp\set{\ol{W} \in M_{n, m-n}(\calO)\big{/}M_{n,m-n}(\calP^{2\ell})}{\check{\pi}^\beta W \equiv 0 \pmod{\calP^{2\ell-1}}} \times N_\ell(\pi^\gamma, \pi^\gamma)\\
& = &q^{2(m-n)\abs{\beta}+2n(m-n)} N_\ell(\pi^\gamma, \pi^\gamma).
\end{eqnarray*}
On the other hand, since the number of $\ol{Y} \in  M_m^{pr}(\calO/\calP^{2\ell}) \cong GL_m(\calO/\calP^{2\ell})$ satisfying \eqref{second cont Y} is equal to $N_\ell^{pr}(\pi^\alp, \pi^\alp) = N_\ell(\pi^\alp, \pi^\alp)$, we see
\begin{eqnarray} 
\mu(\pi^\alp, \pi^\alp) &=&
q^{-\ell m(2m+1)-m(m-1)} N_\ell(\pi^\alp, \pi^\alp) \nonumber\\
&=&
q^{-\ell m(2m+1)-m(m-1) + 2(m-n)\abs{\beta}+2n(m-n)} N_\ell^{pr}(\pi^\beta, \pi^\alp) N_\ell(\pi^\gamma, \pi^\gamma) \nonumber\\
&=&
q^{2(m-n)\abs{\beta}} \cdot \frac{N_\ell^{pr}(\pi^\beta, \pi^\alp)}{q^{\ell n(4m-2n+1) + n(n-1)}} \times \frac{N_\ell(\pi^\gamma, \pi^\gamma)}{q^{\ell (m-n)(2m-2n+1)+ (m-n)(m-n-1)}} \nonumber\\
&=& \label{general form}
q^{2(m-n)\abs{\beta}} \cdot \mu^{pr}(\pi^\beta, \pi^\alp) \cdot \mu(\pi^\gamma, \pi^\gamma).
\end{eqnarray}
Next, assume $\beta_1 < \gamma_{m-n}$. 
For any $V \in M_{
m-n, n}(\calO)$, 
there is $W \in K_n = GL_n(\calO)$ such that $\pi^\beta[W] = \pi^\beta - \pi^{\gamma}[V]$, 
since $\pi^\beta - \pi^{\gamma}[V]$ is $K_n$-equivalent to $\pi^\beta$. 
Then
$$
\begin{pmatrix}V\\W\end{pmatrix}) \in M_{m n}^{pr}(\calO) \; \mbox{and} \; \pi^\alp\left[\begin{pmatrix}V\\W\end{pmatrix}\right] \equiv \pi^\beta \pmod{H_n(\calP, \ell)},
$$
and the number of choice of such $\ol{W} \in M_n(\calO/\calP^{2\ell})$ is equal to $N_\ell(\pi^\beta, \pi^\beta)$. Hence, if $\beta_1 < \gamma_{m-n}$, one has
\begin{eqnarray*}
\mu^{pr}(\pi^\beta, \pi^\alp) &=& q^{-\ell n (4m-2n+1)-n(n-1)} N_\ell^{pr}(\pi^\beta, \pi^\alp) \nonumber\\
&=&
q^{-\ell n (4m-2n+1)-n(n-1)+4\ell n(m-n)} \cdot N_\ell(\pi^\beta, \pi^\beta) \nonumber\\
&=& \mu(\pi^\beta, \pi^\beta),
\end{eqnarray*}
which yields \eqref{eq2 in lem} together with \eqref{general form}.
\qed

\bigskip
By Proposition~\ref{lem: shift} and Proposition~\ref{lem: decomp}, in order to prove Theorem~\ref{th: vol}, it is enough to calculate $\mu(1_n, 1_n)$ and $\mu(h_t, h_t)$, where
\begin{eqnarray} \label{def: ht}
h_t = \pi^{1^{2t}} = \twomatrix{0}{\varPi}{-\varPi}{0} \bot \cdots \bot  \twomatrix{0}{\varPi}{-\varPi}{0} \in X_{2t}.
\end{eqnarray}
We define a $k$-bilinear pairing on the set $\set{X \in M_n(D)}{X^* = X}$ as follows: 
For $B = (b_{ij}),  C = (c_{ij})$, set
\begin{eqnarray}
\pair{B}{C} = \sum_{i=1}^n b_{ii}c_{ii} + \sum_{1 \leq i < j \leq n}\, \Trd(b_{ij}c_{ij}) \in k,
\end{eqnarray}
then we have character sum expressions for $N_\ell^{(pr)}(B, A)$ as follows.

\mslit
\begin{prop}   \label{prop: char sum} 
Let $\ell \geq 1$ and take a character $\chi = \chi_\ell$ of $\calo/\frp^\ell$ such that $\chi$ is nontrivial on $\frp^{\ell-1}/\frp^\ell$.
For $A \in X_m^+$ and $B \in X_n^+$ with $m \geq n$, one has
\begin{eqnarray}  \label{N by character}
N_\ell^{(pr)}(B, A) = q^{-\ell n(2n-1)} 
\sum_{{\scriptsize \begin{array}{c}
\ol{Y} \in M_n(\calO)\big{/}M_n(\calP^{2\ell})\\
Y \equiv Y^*\pmod{\calP^{2\ell}} \end{array} }}  
\sum_{{\scriptsize \begin{array}{c}
\ol{X} \in M_{mn}(\calO)\big{/}M_n(\calP^{2\ell}) \\ X \in M_{mn}^{(pr)}(\calO) \end{array}}} \chi(\pair{A[X]-B}{Y}),
\end{eqnarray}
where $\ol{X}$ and $\ol{Y}$ determine the element $\pair{A[X]-B}{Y}$ in $\calo$ modulo $\frp^{\ell}$. 
\end{prop}

\proof
We write $A[X]-B = (c_{ij})$, and we understand $c_{ij}$'s and entries of $Y$ as elements in $\calO/\calP^{2\ell}$. We calculate the right hand side of the above identity.
\begin{eqnarray*}
\lefteqn{\sum_{{\scriptsize \begin{array}{c}
\ol{Y} \in M_n(\calO/\calP^{2\ell})\\
Y \equiv Y^*\pmod{ \calP^{2\ell}} \end{array} }}  
\sum_{{\scriptsize \begin{array}{c}
\ol{X} \in M_{mn}(\calO)\big{/}M_n(\calP^{2\ell}) \\ X \in M_{mn}^{(pr)}(\calO) \end{array}}} 
\chi(\pair{A[X]-B}{Y})
= \sum_{\ol{X}} \sum_{\ol{Y}} \chi(\pair{(c_{ij})}{Y})} \\
&=&
\sum_{\ol{X}} \prod_{i=1}^n \sum_{ y \in \calo/\frp^\ell}\chi(c_{ii}y) \cdot \prod_{i<j} \sum_{y \in \calO/\calP^{2\ell}} \chi(\Trd(c_{ij}y)) \\
&=&
\sum_{\ol{X}}   \prod_{i=1}^n \left(\begin{array}{ll}
q^{\ell} & \mbox{if } c_{ii} \equiv 0 \pmod{\frp^\ell}\\
0 & \mbox{otherwise}
\end{array} \right) \times 
\prod_{i<j} \left(\begin{array}{ll} 
q^{4\ell} & \mbox{if } c_{ij} \in \calP^{2\ell-1}\\
0 & \mbox{otherwise}
\end{array} \right) \\
&=&
q^{\ell n(2n-1)}N_\ell^{(pr)}(B, A).
\end{eqnarray*}
\qed

For the convenience of later calculation, we note the following.

\begin{prop}   \label{prop: 1n-1n}
\begin{eqnarray}   
&& \label{1-1n}
\mu(1, 1_n) = \mu^{pr}(1,1_n) = 1 - (-q^{-1})^n,\\
&&
\mu(1_n, 1_n) = \prod_{i=1}^n(1-(-q^{-1})^i) = w_n(-q^{-1}). \label{1n-1n}
\end{eqnarray}  
\end{prop}

\proof 
Take $\ell$ to be sufficiently large and $\chi = \chi_\ell$ as in Proposition~\ref{prop: char sum}, For $0 \leq e < \ell$, we set $\chi_{\ell-e}(x) = \chi(\pi^ex)$. Then we may regard $\chi_{\ell-e}$ as a character of $\calo/\frp^{\ell-e}$ that is nontrivial on $\frp^{\ell-e-1}/\frp^{\ell-e}$. 
We may take the representatives of $\calo/\frp^\ell$ as 
\begin{eqnarray} \label{representatives o}
\{0\} \, \bcup\,  \bcup_{e=0}^{\ell-1} \set{\pi^e u}{\ol{u} \in \left(\calo/\frp^{\ell-e}\right)^\times}.
\end{eqnarray}
Then, by \eqref{N by character}, we have
\begin{eqnarray}   
q^\ell N_\ell^{pr}(1, 1_n) 
&=& q^\ell N_\ell(1, 1_n) \nonumber\\
&=&
\sum_{y \in \calo/\frp^\ell}\sum_{{\scriptsize \begin{array}{c}
x_i \in \calO/\calP^{2\ell}\\
1 \leq i \leq n
\end{array} }} \chi ((\sum_{i=1}^n \Nrd(x_i)-1)y) \nonumber \\
&=& \label{for u}
q^{4n\ell} + \sum_{e=0}^{\ell-1} \sum_{\ol{u} \in \left(\calo/\frp^{\ell-e}\right)^\times} 
\left(\sum_{x \in \calO/\calP^{2\ell}} \chi(\pi^eu\Nrd(x))\right)^n \chi(-\pi^eu) \\
&=& 
q^{4n\ell} + \sum_{e=0}^{\ell-1}  \left\{
\left(q^{4e}\sum_{x \in \calO/\calP^{2(\ell-e)}} \chi_{\ell-e}(\Nrd(x))\right)^n \sum_{\ol{u} \in \left(\calo/\frp^{\ell-e}\right)^\times}\chi_{\ell-e}(u) \right\}, \nonumber
\end{eqnarray}
where, since $\Nrd(\calO^\times) = \calo^\times$, one may erase $u$ in the sum with respect to $x$ in \eqref{for u}, and obtain the last expression. 
Since 
\begin{eqnarray}
&&    \label{unit-sum}
\sum_{\ol{u} \in \left(\calo/\frp^m\right)^\times} \chi(u) = \left\{\begin{array}{ll} -1 & \mbox{if } m = 1 \\ 0 & \mbox{if } m > 1 \end{array} \right.,
\end{eqnarray}
we have, as continuation of the above calculation
\begin{eqnarray}
q^\ell N_\ell^{pr}(1, 1_n) 
&=&
q^{4n\ell}-\left(q^{4(\ell-1)}\sum_{x \in \calO/\calP^2} \chi_1(\Nrd(x))\right)^n,
\label{temp1}
\end{eqnarray}
It is easy to see 
\begin{eqnarray}
\sum_{x \in \calO/\calP^2} \chi_1(\Nrd(x)) =  q^2 +  
\sum_{u \in \left(\calo/\frp\right)^\times} \dfrac{q^2(q^2-1)}{q-1} \chi_1(u) 
= 
-q^3, \label{Nrd-sum}
\end{eqnarray}
hence, we obtain by \eqref{temp1} and \eqref{Nrd-sum}
\begin{eqnarray*}
N_\ell^{pr}(1,1_n)&=&
= q^{-\ell}\left(q^{4n\ell}-(-q^{4\ell-1})^n\right)\\
&=&
q^{(4n-1)\ell}(1-(-q^{-1})^n),\\
\mu^{pr}(1, 1_n) &=& (1-(-q^{-1})^n).
\end{eqnarray*}
Finally, by Proposition~\ref{lem: decomp}, we have
\begin{eqnarray*}
\mu(1_n,1_n) = \prod_{r=1}^n \mu^{pr}(1, 1_r) = \prod_{r=1}^n(1-(-q^{-1})^r) = w_n(-q^{-1}).
\end{eqnarray*}
\qed

\bigskip
Next we consider about $\alp = (1, \ldots, 1) \in \Lam_n$, and set $n = 2t$ and $\pi^\alp = h_t$, where
$h_t$ is defined in \eqref{def: ht}. It is convenient to consider the following density
\begin{eqnarray}
N_\ell^{pr}(0, h_t) &=&   \sharp \set{\ol{x} \in M_{n 1}^{pr}(\calO /\calP^{2\ell})}{h_t[x] \equiv 0 \pmod{\frp^\ell}}, \nonumber \\
\mu^{pr}(0, h_t) &=& \lim_{\ell \rightarrow \infty} \frac{N_\ell^{pr}(0, h_t)}{q^{\ell (4n-1)}}.   
\end{eqnarray}
where $M_{n 1}^{pr}(\calO/\calP^{2\ell})$ is defined in \eqref{primitive set}, and in this case
\begin{eqnarray*}
M_{n 1}^{pr}(\calO/\calP^{2\ell}) 
&=&
\set{\ol{x} \in M_{n 1}(\calO/\calP^{2\ell})}{x \notin (\calP)^n}.
\end{eqnarray*}
\begin{lem} Let $n = 2t$. Then    \label{lem: ht-induction}
\begin{eqnarray}
&&
\mu(h_1, h_1) = q^3 \mu^{pr}(0, h_1),\\
&&
\mu(h_t, h_t) = q^{4n-5}\mu^{pr}(0, h_t) \cdot \mu(h_{t-1},h_{t-1}), \quad (t \geq 2).
\end{eqnarray}
\end{lem}

\proof
Take $\ell$ to be sufficiently large. For any $\ol{x} \in M_{n 1}^{pr}(\calO/\calP^{2\ell})$ satisfying $h_t[x] \equiv 0 \pmod{\frp^\ell}$, $x$  can be extended to an element $U \in K_n = GL_n(\calO)$ such that $h_t[U] - h_t \in H_n(\calP, \ell)$. For two such extensions $U$ and $V$ of $x$, we see $h_t[U^{-1}V] - h_t \in H_n(\calP, \ell)$.   Hence the number of such extensions $\ol{U} \in GL_n(\calO/\calP^{2\ell})$ of $\ol{x}$ is equal to the number of $\ol{W} \in GL_n(\calO/\calP^{2\ell})$ such that, when $t \geq 2$, 
\begin{eqnarray}  \label{W in proof}
W  = \left(\begin{array}{cc|c} 1 & b & \; \wt{b} \; \\ 0 & c & \; \wt{c}\; \\
\hline 0 & d & \; D \; \end{array} \right) \in GL_n(\calO), \;  \; h_t[W] \equiv h_t \pmod{H_n(\calP, \ell)},
\end{eqnarray}
where upper left $1, b, 0, c \in \calO$ and other entries are taken with suitable size.  
When $t = 1$, only upper left $2\times2$-block of $W$ in \eqref{W in proof} appears, and we may ignore other entries. We continue the case $t \geq 2$. 
Since 
\begin{eqnarray}   \label{for ht}
h_t[W] = 
\left(\begin{array}{cc|c} 1 & 0 & \; 0 \; \\
b^* & c^* & \; d^* \; \\ \hline \wt{b}^* & \wt{c}^* & D^* \end{array} \right) 
\left(\begin{array}{cc|c} 0 & \varPi c & \varPi \wt{c}\\ -\varPi & -\varPi b & -\varPi \wt{b} \\ \hline 0 & h_{t-1}d & h_{t-1}D \end{array} \right) \equiv  h_t  \pmod{H_n(\calP, \ell)},   \label{htw-1}
\end{eqnarray}
it is easy to see 
\begin{eqnarray} \label{c and c}
c \equiv 1 \pmod{\calP^{2\ell-2}}, \; \wt{c} \equiv 0 \pmod{\calP^{2\ell-2}},
\end{eqnarray}
and the choice of $(c,\wt{c}) \pmod{\calP^{2\ell}}$ is $q^{4(n-1)}$ in $M_{1,n-1}(\calO)$. Then \eqref{for ht} becomes  
\begin{eqnarray}   
h_t[W] \equiv 
\left(\begin{array}{cc|c} 0 & \varPi & \; 0 \; \\ -\varPi & -\Trd(\varPi b)+h_{t-1}[d] & -\varPi\wt{b}+d^*h_{t-1}D\\
\hline
0 & \wt{b}^*\varPi +D^*h_{t-1}d  & h_{t-1}[D] \end{array} \right) \pmod{H_n(\calP, \ell)}, \label{htw-2}
\end{eqnarray}
hence we see
\begin{eqnarray}   
&& \label{D and d}
D \in GL_{n-2}(\calO),  \; h_{t-1}[D] \equiv h_{t-1} \pmod{H_{n-2}(\calP, \ell)}.
\end{eqnarray}
For any $D$ as in \eqref{D and d} and $d \in M_{n-2,1}(\calO)$, we may take $b$ and $\wt{b}$ satisfying $\Trd(\varPi b) \equiv h_{t-1}[d]  \pmod{\frp^\ell}$ and $\varPi \wt{b} \equiv d^*h_{t-1}D \pmod{\calP^{2\ell-1}}$, actually the choice of $b \pmod{\calP^{2\ell}}$ is $q^{3\ell+1}$ in $\calO/\calP^{2\ell}$ and that of $\wt{b} \pmod{\calP^{2\ell}}$  is $q^{4(n-2)}$  in $M_{1,n-2}(\calO/\calP^{2\ell})$. 
If we take $W$ in this way, $\ol{W}$ becomes an element of $GL_n(\calO/\calP^{2\ell})$ since $h_t[W] \equiv h_t \pmod{H_n(\calO, \ell)}$. Hence we see, for $t \geq 2$
\begin{eqnarray}
N_\ell(h_t, h_t) &=& N_\ell^{pr}(0, h_t) \cdot q^{4(n-1)} \cdot N_\ell(h_{t-1}, h_{t-1}) \cdot q^{4\ell(n-2)} \cdot q^{3\ell + 1}\cdot q^{4(n-2)} \nonumber\\
&=&
q^{\ell(4n-5)+8n-11} \cdot N_\ell^{pr}(0, h_t) \cdot N_\ell(h_{t-1},h_{t-1}) \qquad (n = 2t)  \nonumber\\
&=&
q^{\ell(4n-5)+8n-11} \cdot q^{\ell (4n-1)} \mu^{pr}(0, h_t) \cdot q^{\ell (n-2)(2n-3) + (n-2)(n-3)} \mu(h_{t-1}, h_{t-1}) \nonumber\\
&=&
q^{\ell n(2n+1) + n^2+3n-5} \mu^{pr}(0, h_t) \cdot \mu(h_{t-1}, h_{t-1}),
\end{eqnarray}
which yields 
\begin{eqnarray}
\mu(h_t, h_t) = q^{4n-5} \mu^{pr}(0, h_t)\cdot \mu(h_{t-1}, h_{t-1}).
\end{eqnarray}
As for the case $t =1$, we see the condition of $W$ to be $h_1[W] \equiv h_1 \pmod{H_2(\calP, \ell)}$ is $c \equiv 1 \pmod{\calP^{2\ell-2}}$ and $\Trd(\varPi b) \equiv 0 \pmod{\frp^\ell}$, by \eqref{htw-1} and \eqref{htw-2}.  Hence
\begin{eqnarray}
&&
N_\ell(h_1, h_1) = N_\ell^{pr}(0, h_1) \cdot q^4 \cdot q^{3\ell+1},\\ 
&&
\mu(h_1,h_1) = \lim_{\ell \rightarrow \infty} \frac{N_\ell(h_1, h_1)}{q^{10\ell +2}} 
= \lim_{\ell \rightarrow \infty} \frac{N_\ell^{pr}(0, h_1)}{q^{7\ell}} \cdot q^3 \nonumber \\
&&
= q^3 \mu^{pr}(0, h_1).
\end{eqnarray}

\qed

\begin{lem} \label{lem: ht-0}
For each $t \geq 1$, it holds
$$
\mu^{pr}(0, h_t) = q(1-q^{-4t}).
$$
\end{lem}

\proof
We may check the identity \eqref{N by character} holds even when $B = 0$, and we have
\begin{eqnarray}
q^\ell N_\ell(0, h_t) &=& 
\sum_{\ol{z} \in\calo/\frp^\ell} \sum_{\ol{w} \in (\calO/\calP^{2\ell})^{2t}}
\chi_\ell (h_t [w]z) \nonumber \\
&=&
\sum_{\ol{z}\in\calo/\frp^\ell} \left(\sum_{\ol{x}, \ol{y} \in \calO/\calP^{2\ell}} \chi(z \Trd(x\varPi y)) \right)^t \nonumber\\
&=&
q^{1+8t\ell} + \sum_{e=0}^{\ell-2} q^{\ell-e}(1-q^{-1}) \left(\sum_{\ol{x}, \ol{y} \in \calO/\calP^{2\ell}} \chi(\pi^e\Trd(x\varPi y))\right)^t, \label{ht-cal-1}
\end{eqnarray}
where we use the representatives of $\calo/\frp^\ell$ written in \eqref{representatives o} and the property $\Trd$ is $k$-linear. 
To take the summation for $\ol{y}$ in \eqref{ht-cal-1}, we take the representatives of $\calO/\calP^{2\ell}$ for each $e$ with $0 \leq e \leq \ell-2$, as follows:
\begin{eqnarray} \label{decomp of o/pl}
\set{\pi^{\ell-e-1}x}{\ol{x} \in \calO/\calP^{2(e+1)}}\, \bcup\, \bcup_{r=0}^{2\ell-2e-3}\set{\pi^r u}{\ol{u} \in \left(\calO/\calP^{2\ell-r}\right)^\times}.
\end{eqnarray}
Then, we have
\begin{eqnarray}
\lefteqn{\sum_{\ol{x}, \ol{y} \in \calO/\calP^{2\ell}} \chi(\pi^e\Trd(x\varPi y))}\nonumber\\ 
&=&
q^{4(\ell+e+1)} + \sum_{r=0}^{2\ell-2e-3} q^{4\ell-2r}(1-q^{-2}) \sum_{\ol{x}\in \calO/\calP^{2\ell}} \chi(\pi^e\Trd(x\varPi^{r+1})) \nonumber \\
&=&
q^{4(\ell+e+1)}, \label{Trd-cal}
\end{eqnarray}
where we used the fact
\begin{eqnarray}
\sum_{\ol{x}\in \calO/\calP^{2\ell}} \chi(\pi^e\Trd(x\varPi^{r+1})) = 0,
\end{eqnarray}
which holds since $2e+r+1 < 2\ell-1$ and $\chi$ is nontrivial on $\frp^{\ell-1}/\frp^\ell$.    
Hence we obtain
\begin{eqnarray}
q^\ell N_\ell(0, h_t) &=& 
q^{1+8t\ell} + \sum_{e=0}^{\ell-2} q^{\ell-e}(1-q^{-1})q^{4t(\ell+e+1)}.   \label{ht-0-1}
\end{eqnarray}

Next, we calculate the number $N_\ell^{imp}(0,h_t)$ of imprimitive solutions for $h_t[x] \equiv 0 \pmod{\frp^\ell}$ by character sum as follows
\begin{eqnarray} 
q^\ell N_\ell^{imp}(0, h_t) &=& 
q^\ell \sharp \set{\ol{x} \in M_{n,1}(\calP/\calP^{2\ell})}{h_t[x] \equiv 0 \pmod{\frp^\ell}} \nonumber\\
&=&
\sum_{\ol{z}\in\calo/\frp^\ell} \sum_{\ol{w} \in (\calP/\calP^{2\ell})^{2t}}
\chi_\ell (h_t [w]z) \nonumber \\
&=&
\sum_{\ol{z}\in\calo/\frp^\ell} \left(\sum_{\ol{x}, \ol{y} \in \calP/\calP^{2\ell}} \chi(z \Trd(x\varPi y)) \right)^t \nonumber \\
&=&
q^{2+4t(2\ell-1)} + \sum_{e=0}^{\ell-3} q^{\ell-e}(1-q^{-1}) \left(\sum_{\ol{x}, \ol{y} \in \calP/\calP^{2\ell}} \chi(\pi^e\Trd(x\varPi y))\right)^t, \label{ht-cal-im-1}
\end{eqnarray}
where we use \eqref{representatives o} and the fact $\Nrd(x\varPi y) \in \frp^2$ for $x, y \in \calP$.
In the similar way to calculate \eqref{Trd-cal}, we have for each $e$ with $0 \leq e \leq \ell-3$,
\begin{eqnarray}
\lefteqn{\sum_{\ol{x}, \ol{y} \in \calP/\calP^{2\ell}} \chi(\pi^e\Trd(x\varPi y))} \nonumber\\
&=&
q^{4\ell-2+4e+6} + \sum_{r=1}^{2\ell-2e-4} q^{4\ell-2r}(1-q^{-2}) \sum_{\ol{x}\in \calP/\calP^{2\ell}} \chi(\pi^e\Trd(x\varPi^{r+1})) \nonumber\\
&=&
q^{4(\ell+e+1)}+ \sum_{r=1}^{2\ell-2e-4} q^{4\ell-2r}(1-q^{-2}) \cdot 0 \nonumber\\ 
&=&
q^{4(\ell+e+1)}. \label{ht-cal-im-2}
\end{eqnarray}
Hence we have 
\begin{eqnarray}
q^\ell N_\ell^{imp}(0, h_t) &=& 
q^{2+4t(2\ell-1)} + \sum_{e = 0}^{\ell-3} q^{\ell-e}(1-q^{-1})q^{4t(\ell+e+1)}   \label{ht-0-2}
\end{eqnarray}
By \eqref{ht-0-1} and \eqref{ht-0-2}, we obtain
\begin{eqnarray}
N_\ell^{pr}(0, h_t) &=& 
q^{(8t-1)\ell}q(1-q^{-4t}),\quad (n=2t)
\end{eqnarray}
which yields $\mu^{pr}(0, h_t) =  q(1-q^{-4t})$.   
\qed

\bigskip
By Lemma~\ref{lem: ht-induction} and Lemma~\ref{lem: ht-0}, we have the following.

\begin{prop} \label{prop: ht-ht}
Let $n = 2t$ and $\alp = (1, \ldots, 1) \in \Lam_n^+$.  Then $\pi^\alp = h_t$, and  it holds
\begin{eqnarray}   \label{ht-ht formula}
\mu(h_t, h_t) = q^{4t^2} \cdot \prod_{i=1}^{t} (1-q^{-4i}) = q^{n^2}\cdot w_{\frac{n}{2}}(q^{-4}). 
\end{eqnarray}
\end{prop}

\medskip
{\it Proof of} Theorem~\ref{th: vol}.\; 
Take $\alp \in \Lam_n^+$ as in \eqref{shape of alp}.  Then, by Propositions~\ref{lem: shift}, \ref{lem: decomp}, \ref{prop: 1n-1n}, and \ref{prop: ht-ht}, we see
\begin{eqnarray*}
\mu(\pi^\alp, \pi^\alp) &=& 
q^{\sum_{j=2}^r 2(m_1+\cdots +m_{j-1})m_j\gamma_j} \cdot \prod_{j=1}^r\mu(\pi^{\gamma_j^{m_j}}, \pi^{\gamma_j^{m_j}}) \nonumber\\ 
%
&=&
q^{m_\alp} \cdot 
\prod_{{\scriptsize \begin{array}{c} 1 \leq j \leq r\\ 2\mid \gamma_j \end{array} }} w_{m_j}(-q^{-1}) 
\cdot
\prod_{{\scriptsize \begin{array}{c} 1 \leq j \leq r\\ 2\not\vert \gamma_j \end{array} }} w_{\frac{m_j}{2}}(q^{-4}),
\end{eqnarray*}
where
\begin{eqnarray*}
m_\alp &=& \sum_{i=2}^r 2(m_1+\cdots +m_{j-1})m_j\gamma_j + 
\sum_{{\scriptsize \begin{array}{c} 1 \leq j \leq r\\ 2\mid \gamma_j \end{array} }} \frac{\gamma_j}{2}m_j(2m_j-1) \\
&& \qquad+ 
\sum_{{\scriptsize \begin{array}{c} 1 \leq j \leq r\\ 2\not\vert \gamma_j \end{array} }} (\frac{\gamma_j-1}{2} m_j(2m_j-1) + m_j^2)\\
&=&
2\sum_{j=2}^r(m_1+\cdots +m_{j-1})m_j\gamma_j + \sum_{j=1}^r(m_j(m_j-1) +\frac12m_j)\gamma_j + 
\frac12 \sum_{{\scriptsize \begin{array}{c}  1 \leq i \leq r \\ 2\not\vert \gamma_i \end{array} }} m_j \\
&=&
2n(\alp) + \frac12\abs{\alp} + \frac12\sharp\set{i}{\alp_i \mbox{ is odd}}. 
\end{eqnarray*}
For $\wt{\alp}=\alp+(2e,\ldots,2e) \in \Lam_n$, we have $\mu(\pi^{\wt{\alp}},\mu^{\wt{\alp}}) = q^{en(2n-1)}\mu(\pi^\alp, \mu^\alp)$ by  Remark~\ref{rem: shift}. On the other hand, we have $2n(\wt{\alp})+\frac12\abs{\wt{\alp}} = 2n(\alp)+\frac12 \abs{\alp} + en(2n-1)$, hence we see \eqref{density oneself} holds for any $\alp \in \Lam_n$
\qed

\slit
{\bf 2.3.} We introduce the Schwartz space $\SKX$ by 
$$
\SKX = \set{\vphi: X \longrightarrow \C}{\mbox{left $K$-invariant, compactly supported}},
$$ 
that is spanned by the characteristic functions of $K\cdot \pi^\alp, \; \alp \in \Lam_n$ over $\C$. 
It is an $\hec$-submodule of $\CKX$ (cf. \eqref{CKX}, \eqref{hec acts on CKX}).
We define an integral transform $F_0$ on $\SKX$ as follows:
\begin{eqnarray}  \label{sph trans F0}
\begin{array}{lccl}
F_0 : & \SKX & \longrightarrow & \C(q^{s_1}, \ldots, q^{s_n}) , \\[2mm]
& \vphi & \longmapsto & \int_{X}\vphi(x) \omega(x^{-1}; s) dx, 
\end{array}
\end{eqnarray}
where $dx$ is a $G$-invariant measure on $X$. 
We call $F_0$ a spherical Fourier transform on $\SKX$, and we will normalize $F_0$ suitably and define the spherical transform $F$ in \S 3.

\begin{prop} \label{prop: sph tr}
The spherical Fourier transform $F_0$ defined in \eqref{sph trans F0} is injective and compatible with the action of $\hec$: 
$$
F_0(f*\vphi)(s) = \lam_s(f)F_0(\vphi), \quad f \in \hec, \; \vphi \in \SKX,
$$
where $\lam_s$ is defined in \eqref{lam_s}.
\end{prop}

\medskip
The injectivity of $F_0$ is proved in a similar way to the cases of other sesquilinear forms by using Lemma~\ref{lem: key} below and Theorem~\ref{th: ind} (cf. \cite[\S 3 Theorem]{H1JapJM}). 
We define a binary relation $\succ$ on $\Lam_n$ by 
\begin{eqnarray*}
\lam \succ \mu \, \Llra\,  \begin{array}{l}
\lam = \mu, \mbox{or there is some $t$ with $1 \leq t \leq n-1$ satisfying }\\
\qquad  \lam_{n-t} > \mu_{n-t}, \; \lam_{n-t+1}=\mu_{n-t+1}, \ldots, \lam_n=\mu_n.
\end{array} 
%
\end{eqnarray*}

\begin{lem} \label{lem: key}
Let $n \geq 2$. For any $\alp \in \Lam_n^+$, there exists $\beta \in \Lam_{n-1}^+$ such that\\
\qquad{\rm (i)} $\mu^{pr}(\pi^\beta, \pi^\alp) \ne 0$, and \\
\qquad{\rm (ii)} if $\gamma \in \Lam_n^+$ satisfies 
$\abs{\gamma} = \abs{\alp}, \; \gamma \succ \alp$ and $\mu^{pr}(\pi^\beta, \pi^\gamma) \ne 0$,
then $\gamma = \alp$.
\end{lem} 

Similar lemma was introduced first by Kitaoka (\cite{Kitaoka}) for symmetric forms and by the author for hermitian forms (\cite[\S3]{H1JapJM}), and the above lemma can be proved similarly, so we note here that one may take $\beta \in \Lam_{n-1}$ as 
$\beta = (\alp_2, \alp_3, \ldots, \alp_n)$ (resp. $(\alp_2+1, \alp_3, \ldots, \alp_n)$ if $\alp_1$ is even (resp. odd).

\bigskip
{\it Proof of} Proposition~\ref{prop: sph tr}.\;
Let $f \in \hec$ and $\vphi \in \SKX$. Then we have
\begin{eqnarray*}
F_0(f*\vphi)(s) &=&  \int_X \int_G f(g)\vphi(g^{-1}\cdot x)\omega(x^{-1};s)dgdx\\
&=&
\int_G \int_X f(g)\vphi(y)\omega(g^{* -1}\cdot y^{-1}; s)dydg\quad (y = g^{-1}\cdot x)\\
&=&
\int_X \int_G f(g^*)\omega(g^{* -1}\cdot y^{-1}; s) dg \vphi(y)dy \\
&=&
\lam_s(f) \int_X \vphi(y)\omega(y^{-1}; s) = \lam_s(f) F_0(\vphi).
\end{eqnarray*}
We prove the injectivity of $F_0$ by induction on $n$. For $\alp \in \Lam_n$, we denote by $\vphi_\alp \in \SKX$ the characteristic function of $K \cdot (\pi^\alp)^{-1}$. Then we have $F_0(\vphi_\alp) = volume(K\cdot (\pi^\alp)^{-1}) \times \omega(\pi^\alp; s) \ne 0$,  hence the injectivity of $F_0$ is equivalent to the linear independence of $\omega(\pi^\alp; s)$ for $\alp \in \Lam_n$.  
It is clear that $F_0$ is injective for $n=1$. We assume that $F_0$ is injective for $n-1$ and not injective for $n$, and take $0 \ne \vphi \in {\rm Ker}(F_0)$. 
We may assume $\vphi$ is of the following shape:
\begin{eqnarray}  \label{shape of vphi}
\vphi = \sum_{i=1}^\ell\, c_i \vphi_{\alp^{(i)}}, \quad \ell \geq 2,  \; c_i \ne 0, \; \alp^{(i)} \in \Lam_n^+ \; (1 \leq i \leq \ell), \; \alp^{(i)} \ne \alp^{(j)} \mbox{ if } i\ne j.
\end{eqnarray}
Since we have, for any $\alp \in \Lam_n$, 
\begin{eqnarray*}
\omega(\pi^\alp; s) = q^{\frac{\abs{\alp}}{2}s_n} \omega(\pi^{\alp}; s_1, \ldots, s_{n-1},0), \quad \omega(\pi^{\alp}; s_1, \ldots, s_{n-1},0) \in \C(q^{s_1},\ldots,q^{s_{n-1}}),
\end{eqnarray*}
looking at the exponent of $q^{s_n}$ in $\omega(\pi^\alp; s)$, we may assume that $\abs{\alp^{(i)}} = \abs{\alp^{(1)}}$ for any $i$. Assume that $\alp^{(1)}$ is the smallest with respect to the order $\succ$ within $\set{\alp^{(i)}}{1 \leq i \leq \ell}$. Since $F_0(\vphi) = 0$, we obtain by Theorem~\ref{th: ind}
\begin{eqnarray} \label{n-1-rel}
\sum_{i=1}^\ell c_i d_i \sum_{\gamma \in \Lam_{n-1}^+}\, \frac{\mu^{pr}(\pi^\gamma, \pi^{\alp^{(i)}})}{\mu(\pi^\gamma, \pi^\gamma)} \omega(\pi^\gamma; s_1, \ldots, s_{n-1}) = 0,
\end{eqnarray}
where $d_i (> 0)$ is the volume of $K\cdot(\pi^{\alp^{(i)}})^{-1}$. By induction hypothesis that $F_0$ is injective for $n-1$, we see by \eqref{n-1-rel}
\begin{eqnarray}
\sum_{i=1}^\ell c_id_i \mu(\pi^\gamma, \pi^{\alp^{(i)}}) = 0, \quad \mbox{for every }\gamma \in \Lam_{n-1}^+.
\end{eqnarray}
For $\beta \in \Lam_{n-1}^+$ associated with $\alp^{(1)}$ in Lemma~\ref{lem: key}, we see that
$\mu(\pi^\beta, \pi^{\alp^{(1)}}) \ne 0$ and $\mu(\pi^\beta, \pi^{\alp^{(i)}}) = 0$ for $i \ne 1$, and we obtain $c_1=0$, which contradicts \eqref{shape of vphi}. \qed

\vspace{1cm}
\newpage
\Section{Functional equations of spherical functions} 

First we note the result for size $2$, which follows from Theorem~\ref{th: ind} with some calculation of local densities (cf. \cite{OY}).

\begin{prop} \label{prop: size2}
For any $\alp \in \Lam_2$, one has 
\begin{eqnarray*}
\lefteqn{\omega(\pi^\alp; z)} \nonumber\\
& = &
\left\{\begin{array}{ll}\
\displaystyle{\frac{q^{\pair{\lam}{z_0}}}{1+q^{-2}}\cdot \frac{1}{q^{z_2}-q^{z_1+1} } \sum_{\sigma \in S_2}}\sigma\left(q^{\pair{\lam}{z}}\frac{(q^{z_1}- q^{z_2-2})(q^{z_1}-q^{z_2+1})}{q^{z_1}-q^{z_2}} \right) & \mbox{if } \; \alp = 2\lam,\\
q(1-q^{-1})\dfrac{q^{e(z_1+z_2)}}{q^{z_2}-q^{z_1+1}} & \mbox{if } \; \alp = (2e-1, 2e-1),
\end{array} \right.
\end{eqnarray*}
where $z_0 = (-1, 1)$ corresponds to the $s$-variable ${\bf0} = (0,0)$, $\pair{\lam}{z} = \lam_1z_1 + \lam_2z_2$ and $S_2$ acts on $\{z_1, z_2\}$ by permutation. 
Especially, for any $x \in X_2$, one has 
\begin{eqnarray} \label{sym2}
(q^{z_2}-q^{z_1+1})\cdot \omega(x; z) \in \C[q^{\pm z_1}, q^{\pm z_2}]^{S_2}.
\end{eqnarray}
\end{prop}

\bigskip
The property \eqref{sym2} follows from the explicit formula for $\omega(\pi^\alp; z)$,
since any $x \in X_2$ belongs to some orbit $K_2\cdot \pi^\alp, \; \alp \in \Lam_2$ and 
$\omega(x; z) = \omega(\pi^\alp; z)$.

\bigskip
For the study of the functional equations and holomorphy of $\omega(x; s)$ for general $n$,  we use the same strategy used in the case of unramified hermitian forms. We introduce the following integral for  $\xi \in \SKX$ 
\begin{eqnarray}
\Phi(s, \xi)  
= \int_{X} \abs{\bfd(x)}^s \xi(x) dx, \quad
\abs{\bfd(x)}^s = \left\{\begin{array}{ll}
\prod_{i=1}^n \abs{d_i(x)}^{s_i} & \mbox{if} \; x \in X^{op}\\
0 & \mbox{otherwise},
\end{array} \right.
\end{eqnarray}
where $dx$ is a $G$-invariant measure on $X$. The above integral is a finite linear sum of spherical functions $\omega(x; s)$, hence it is absolutely convergent if $\real(s_i) \geq 0, \; 1 \leq i \leq n-1$, and continued to a rational function of $q^{s_1}, \ldots, q^{s_n}$.  Keeping the relation \eqref{change of var} between $s$ and $z$, we denote $\Phi(z, \xi)$.

\begin{lem}  \label{lem: for funeq}
Let $n \geq 2$ and take $\alp$ with $1 \leq \alp \leq n-1$. Assume that $\real(s_i) \geq 0, \; 1 \leq i \leq n-1$. Then for any $\xi \in \SKX$, the following identity holds 
\begin{eqnarray}   \label{fund eq for funeq}
\Phi(s, \xi) = 
\int_{X^{op}} \prod_{i \ne \alp, \alp+1} \abs{d_i(x)}^{s_i} \cdot \prod_{j=\alp\pm 1}\abs{d_j(x)}^{\frac{s_\alp}{2}+s_j} \cdot \xi(x) \cdot \omega^{(2)}(\wt{x}; s_\alp, -\frac{s_\alp}{2}) dx,
\end{eqnarray}
where $\wt{x}$ is the lower right $(2\times2)$-block of $(x^{(\alp+1)})^{-1}$ and $\omega^{(2)}(y;s)$ indicates the spherical function of size $2$.  
\end{lem}

\proof
Take any $\alp$ with $1 \leq \alp \leq n-1$ and $\xi \in \SKX$. We assume that $\real(s_i) \geq 0, \; 1 \leq i \leq n-1$.
We define an embedding $\iota = \iota_\alp$ from $K_2 = GL_2(\calO)$ into $K = K_n$ by
\begin{eqnarray*}
\iota: K_2 \longrightarrow K, \; k \longmapsto \left(\begin{array}{c|c|c}
1_{\alp-1}&0&0\\
\hline
0&k&0\\
\hline
0&0&1_{n-\alp}
\end{array}\right),
\end{eqnarray*}
and consider the integral
\begin{eqnarray}
\Phi(s,\xi)&=&
\int_{X} \abs{\bfd(x)}^s\int_{K_2}\xi(\iota(k)^{-1}\cdot x)dkdx \nonumber\\
&=&
\int_{K_2}\int_{X} \abs{\bfd(\iota(k)\cdot x)}^s \xi(x)dxdk. \label{first step}
\end{eqnarray}
Here we recall a well known fact on miner determinants of matrices over a field $F$: 
\begin{eqnarray} \label{miner det}
\det(A^{(i)}) = \det(A)\det({A^{-1}}_{(n-i)}), \quad \mbox{for \;} A \in GL_n(F),
\end{eqnarray}
where $A^{(i)}$(resp. $A_{(j)}$) indicates the upper left $(i \times i)$-block (resp. the lower right $(j\times j)$-block) of $A$.
In the present case, we consider the reduced norm $\Nrd: D \longrightarrow k$ and relative invariants $d_i(x)$, which satisfy $d_i(x)^2 = \Nrd(x^{(i)})$ on $X$, as introduced in \S1.  Hence we have
\begin{eqnarray*}
d_i(\iota(k)\cdot x) &=& d_i(x) \; \mbox{unless \; } i = \alp,\\
\Nrd((\iota(k)\cdot x)^{(\alp)}) & =& \Nrd(x^{\alp+1}) \Nrd( {\left( (\iota(k)\cdot x)^{(\alp+1)}\right)^{-1}}_{(1)}),\\
d_\alp(\iota(k)\cdot x) &=& d_{\alp+1}(x) d_1(j_2k^{* -1}\cdot \wt{x}).
\end{eqnarray*}
We continue the calculation \eqref{first step} as follows
\begin{eqnarray*}
\Phi(s,\xi)&=&
\int_{X^{op}} \xi(x) \prod_{i \ne \alp, \alp+1}\abs{d_i(x)}^{s_i} \cdot \abs{d_{\alp+1}(x)}^{s_\alp+s_{\alp+1}} \int_{K_2}\abs{d_1(j_2k^{*-1}\cdot \wt{x})}^{s_\alp}dkdx\\
&=&
\int_{X^{op}} \xi(x) \prod_{i \ne \alp, \alp+1}\abs{d_i(x)}^{s_i} \cdot \abs{d_{\alp+1}(x)}^{s_\alp+s_{\alp+1}} \cdot \omega^{(2)}(\wt{x}; s_\alp, 0)dx.
\end{eqnarray*}
By the definition of $\wt{x}$, we see
\begin{eqnarray*}
d_2(\wt{x}) = d_{\alp-1}(x)d_{\alp+1}(x)^{-1},
\end{eqnarray*}
and we obtain
\begin{eqnarray*}
\Phi(s,\xi)&=&
\int_{X^{op}} \xi(x) \prod_{i \ne \alp, \alp\pm1}\, \abs{d_i(s)}^{s_i} \cdot \prod_{j=\alp\pm1}\, \abs{d_j(x)}^{s_j+\frac{s_\alp}{2}} \cdot \omega^{(2)}(\wt{x}; s_\alp, -\frac{s_\alp}{2})dx.
\end{eqnarray*}
\qed

%
\begin{prop}
Under the relation \eqref{change of var} of $s$ and $z$, the function 
\begin{eqnarray*}
\prod_{1 \leq i < j \leq n} (q^{z_j}-q^{z_i+1}) \times \Phi(s, \xi), \quad (\xi \in \SKX)
\end{eqnarray*}
is holomorphic in $\C^n$ and $S_n$-invariant, hence it is an element of 
\begin{eqnarray*}
\C[q^{\pm z_1}, \ldots, q^{\pm z_n}]^{S_n}.
\end{eqnarray*}
\end{prop}

\proof
By the relation \eqref{change of var}, $S_n$ acts on variable $s$ as follows. Let $\sigma_\alp = (\alp \; \alp+1) \in S_n, \; 1 \leq \alp \leq n-1$, then its action on $\{s_i\}$ is
\begin{eqnarray}
&&
\sigma_\alp(s_i) = s_i, \quad \mbox{unless }\; i = \alp, \alp\pm1,\nonumber\\
&&
\sigma_\alp(s_\alp) = -s_\alp-4, \nonumber\\
&&
\sigma(s_j) = s_j+s_\alp+2 \quad \mbox{if }\; j = \alp\pm1.\label{sigma-alp}
\end{eqnarray}
Set $\calD_\alp=\calD_0 \cup \calD_{\alp,1}\cup \calD_{\alp,2}$, where
\begin{eqnarray}
&&
\calD_0 = \set{s \in \C^n}{\real(s_i) \geq 0, \; (1\leq i \leq n-1)},\nonumber \\[2mm]
&&
\calD_{\alp,1}=\sigma_\alp(\calD_0)= \set{s \in \C^n}{\begin{array}{l}
\real(s_i)\geq 0 \; \mbox{for } i \in [1, n-1], i \ne \alp, \alp\pm 1\\
\real(s_\alp) \leq -4\\
\real(s_\alp+s_j+2) \geq 0 \; \mbox{for } j = \alp\pm 1 \in [1, n-1]
\end{array}},\nonumber \\[2mm]
&&
\calD_{\alp,2}=\set{s \in \C^n}{\begin{array}{l}
\real(s_i) \geq 0, \; \mbox{for } i \in [1, n-1], i \ne \alp, \alp\pm 1,\\
-4 \leq \real(s_\alp) \leq 0,\\
\real(s_\alp/2+s_j) \geq 0, \; \mbox{for } j = \alp\pm 1 \in [1, n-1]
\end{array}}. \label{calD}
\end{eqnarray}
By the relation of $s$ and $z$ and Proposition~\ref{prop: size2}, one has
\begin{eqnarray*}
&&
q^{-\frac{z_\alp + z_{\alp+1}}{2}}(q^{z_{\alp+1}}-q^{z_\alp+1}) = q^{\frac{s_\alp}{2}+1}-q^{-\frac{s_\alp}{2}},\\
&&
(q^{\frac{s_\alp}{2}+1}-q^{-\frac{s_\alp}{2}})\omega^{(2)}(x; s_\alp, -\frac{s_\alp}{2}) = 
(q^{-\frac{s_\alp}{2}-1}-q^{\frac{s_\alp}{2}+2})\omega^{(2)}(x; -s_\alp-4, \frac{s_\alp}{2}+2) \in \C[q^{\pm\frac{s_\alp}{2}}]^{\gen{\sigma_\alp}}.\nonumber\\
\label{first eq}
\end{eqnarray*}
Then by Lemma~\ref{lem: for funeq}, one has for $s \in \calD_0$
\begin{eqnarray}
(q^{z_{\alp+1}}-q^{z_\alp+1})\Phi(s, \xi) &=& q^{\frac{z_\alp + z_{\alp+1}}{2}}
\int_{X^{op}} \prod_{i \ne \alp, \alp+1} \abs{d_i(x)}^{s_i} \cdot \prod_{j=\alp\pm 1}\abs{d_j(x)}^{\frac{s_\alp}{2}+s_j} \cdot \xi(x) \nonumber \\
&& \label{for feq2}
\quad \times (q^{\frac{s_\alp}{2}+1}-q^{-\frac{s_\alp}{2}})\omega^{(2)}(\wt{x}; s_\alp, -\frac{s_\alp}{2}) dx.
\end{eqnarray}
Since the integrand of RHS of \eqref{for feq2} is $\sigma_\alp$-invariant, we see the above integral is absolutely convergent for $s \in \sigma_\alp(\calD) = \calD_{\alp, 1}$.  
The region $\calD_{\alp,2}$ is $\sigma_\alp$ -invariant, and we see that
\begin{eqnarray*} 
&&
\prod_{i \ne \alp, \alp+1} \abs{d_i(x)}^{s_i} \cdot \prod_{j=\alp\pm 1}\abs{d_j(x)}^{\frac{s_\alp}{2}+s_j} \mbox{ is bounded for }\; s \in \calD_{\alp, 2},\\
&&
(q^{\frac{s_\alp}{2}+1}-q^{-\frac{s_\alp}{2}})\omega^{(2)}(x; s_\alp, -\frac{s_\alp}{2}) \mbox{  is a polynomial in } \; q^{\pm\frac{s_\alp}{2}}.
\end{eqnarray*}
Since $\xi$ is compactly supported, RHS of \eqref{for feq2} is absolutely convergent also for $s \in \calD_{\alp,2}$, and so $(q^{z_{\alp+1}}-q^{z_\alp+1})\Phi(s, \xi)$ is holomorphic in $\calD_\alp$ and $\sigma_\alp$-invariant.
Since  
\begin{eqnarray*}
\prod_{1\leq i<j \leq n}\, (q^{z_j}-q^{z_i}) \big{/} (q^{z_{\alp+1}}-q^{z_\alp+1}) 
\end{eqnarray*}
is $\sigma_\alp$-invariant for ant $\alp$ and holomorphic in $\C^n$,
\begin{eqnarray*}
\prod_{1\leq i<j \leq n}(q^{z_j}-q^{z_i}) \times \Phi(s, \xi) \mbox{ is holomorphic in }
\calC = \bcup_{\alp=1}^{n-1} \calD_\alp \mbox{ and $S_n$-invariant}.
\end{eqnarray*}
Hence 
\begin{eqnarray} \label{normalized Phi(s)}
\prod_{1\leq i<j \leq n}(q^{z_j}-q^{z_i}) \times \Phi(s, \xi)
\end{eqnarray}
is holomorphic in 
$\bcup_{\sigma \in S_n} \sigma(\calC)$ and its convex hull $\C^n$ and $S_n$-invariant. Since \eqref{normalized Phi(s)} is a rational function of $q^{z_1},\ldots, q^{z_n}$, we see that it is a symmetric Laurent polynomial, thus we have
\begin{eqnarray*}
\prod_{1\leq i<j \leq n}(q^{z_j}-q^{z_i}) \times \Phi(s, \xi) \in \C[q^{\pm z_1},\ldots,  q^{\pm z_n}]^{S_n}.
\end{eqnarray*}
\qed

\bigskip
Taking the characteristic function of $K\cdot x$ for $x \in X_n$ as $\xi$, we obtain the following theorem. 

\bigskip
\begin{thm}   \label{th: feq}
Set 
\begin{eqnarray} \label{def G_n(z)}
\Psi(x;z) = G_n(z) \cdot \omega(x; z), \quad G_n(z) = \prod_{1 \leq i < j \leq n} (q^{z_j}-q^{z_i+1}), 
\end{eqnarray}
then $\Psi(x;z)$ is holomorphic and $S_n$-invariant spherical function on $X$, thus 
\begin{eqnarray*}   \label{def Psi(x;z)}
\Psi(x;z) \in \C[q^{\pm z_1}, \ldots, q^{\pm z_n}]^{S_n}.
\end{eqnarray*}
\end{thm}

\bigskip
In consideration of Theorem~\ref{th: feq}, we normalize the spherical Fourier transform $F_0$ given in \eqref{sph trans F0} as follows:
\begin{eqnarray}  \label{sph transf-s}
\begin{array}{cccl}
F : & \SKX & \longrightarrow &\C[q^{\pm z_1}, \ldots, q^{\pm z_n}]^{S_n} (= \calR, \mbox{say})\\
& \vphi & \longmapsto & \widehat{\vphi}(z) = \displaystyle{\int_{X}\vphi(x) \cdot \Psi(x^{-1}; z) dx}.
\end{array} 
\end{eqnarray}
Then we obtain the following theorem by Theorem~\ref{th: feq}.

\begin{thm} \label{th: spTrans}
The normalized spherical Fourier transform $F$ is an injective $\hec$-module map, hence one has the commutative diagram
\begin{eqnarray}   \label{comm diag}
\begin{array}{ccccc}
\hec &\times & \SKX & \stackrel{*}{\longrightarrow} & \SKX\\
\mapdownlr{\lam_z\,}{\rotatebox{90}{$\sim$}} && \mapdownl{F\;} & \circlearrowleft & \mapdownl{F\;} \\
\calR & \times &\calR& \longrightarrow & \calR,
\end{array}
\end{eqnarray}
where the upper $*$ is the action of $\hec$ on $\SKX$, the lower arrow is the multiplication in $\calR$, and $\lam_z$ is the Satake isomorphism defined in \eqref{lam_z}.
\end{thm}

\vspace{1cm}
\Section{Explicit formula for $\omega(x;z)$}
{\bf 4.1.} As for the explicit formula of $\omega(x; z)$, it suffices to determine it at each representatives of $K$-orbit in $X$, i.e. at each $\pi^\alp, \;\alp \in \Lam_n$(cf. \eqref{Lam-n}). 
We may apply the general expression formula of spherical function on homogeneous spaces (cf. \cite[Prop.1.9]{JMSJ}, \cite[\S 2]{French}). In the present case the situation becomes simpler, since $\omega(x;z)$ has good functional equations and $X^{op}$ is a single $B$-orbit, and all the assumptions to apply the general expression formula of $\omega(x;s)$ are satisfied: i.e. (i) $X$ has only a finite numbers of $B$-orbits over the algebraic closure of $k$; (ii)  the relative invariants $d_i(x), \; 1 \leq i \leq n$ are
regular functions on $X$ and corresponding characters $\psi_i$ generate the group $\frX(B)$ of rational characters of $B$ defined over $k$; (iii) for $y \in X \backslash X^{op}$, there exists a rational $\psi \in \frX(B)$ whose restriction to the identity component of the stabilizer $B_y$ is not trivial.

\mslit
For each $\alp = (\alp_i) \in \Lam_n$, we set
\begin{eqnarray} \label{def: lam-alp}
\lam_\alp = (\lam_i) \in \wt{\Lam}_n \quad \mbox{by}\quad \lam_i = \left[\dfrac{\alp_i+1}{2}\right],
\end{eqnarray}
where $[\;\;]$ is the Gauss symbol.
If $\alp$ has an odd entry, odd entries appear in pairs. We assume they are 
\begin{eqnarray}
\alp_{\ell_1}, \alp_{\ell_1+1}, \ldots, \alp_{\ell_k}, \alp_{\ell_k+1}, \quad \ell_1<\ell_2< \cdots < \ell_k,
\end{eqnarray}
and set
\begin{eqnarray} \label{I_odd}
I_{odd}(\alp) = \{\ell_1, \ldots, \ell_k\}, \quad  c_{odd}(\alp) = (1-q^{-1})^k \cdot q^{\sum_{\ell \in I_{odd}(\alp)}\, (n-2\ell+1)}.
\end{eqnarray}
If  $\alp$ has no odd entry we say $\alp$ is {\it even}, and set $I_{odd}(\alp) = \emptyset$ and $c_{odd}(\alp) = 1$ for convenience. Only if $\alp$ is even, $\pi^\alp$ is diagonal and $\lam_\alp = \dfrac{\alp}{2}$.  
We define a paring on $\Z^n \times \C^n$ as follows:
\begin{eqnarray}
\pair{\lam}{z} = \sum_{i=1}^n\, \lam_iz_i, \quad (\lam \in \Z^n, \; z \in \C^n).
\end{eqnarray}


\medskip
\begin{thm} $($Explicit Formula$)$\label{th: explicit}
For any $\alp \in \Lam_n$, one has
\begin{eqnarray}
&&
\Psi(\pi^\alp; z) = \omega(\pi^\alp; z) \cdot G_n(z) \\
&& 
=\frac{(1-q^{-2})^n \cdot c_{odd}(\alp) \cdot q^{\pair{\lam_\alp}{z_0}}} {w_n(q^{-2}) }
\sum_{\sigma \in S_n}\, 
\sigma \left( \frac{ q^{\pair{\lam_\alp}{z}} }{ \displaystyle{\prod_{\ell \in I_{odd}(\alp)}} (q^{z_{\ell}} - q^{z_{\ell+1}+1 })} \prod_{i<j} \frac{(q^{z_i}-q^{z_j+1})(q^{z_i}-q^{z_j-2}) }{ q^{z_i}-q^{z_j}} \right), \nonumber
\end{eqnarray}
where $w_n(t) = \prod_{i=1}^n(1-t^i)$, 
$
G_n(z) = \prod_{1 \leq i < j \leq n}\, (q^{z_j}-q^{z_i+1})$ $($given in Theorem~\ref{th: feq}$)$ and $z_0 = (-n+1, -n+3, \ldots, n-1) \in \C^n$ is the corresponding value in $z$-variable to $s = \bf0 \in \C^n$.
\end{thm} 

\proof
Applying \cite[Prop.1.9]{JMSJ} to the present case, we obtain for generic $z$, 
\begin{eqnarray} \label{4-first step}
\omega(x;z) &=& \frac{1}{Q_n} \sum_{\sigma \in S_n}\, \gamma(\sigma(z)) \Gamma_\sigma(z) \int_{U}\abs{\bfd(\nu\cdot x)}^{\sigma(s)}d\nu,
\end{eqnarray}
where $U$ is the Iwahori subgroup associated with $B$, $d\nu$ is the Haar measure on $U$, $Q_n$ and $\gamma(z)$ are determined by the group $GL_n(D)$ as follows, and $\Gamma_\sigma(z)$ is determined by the functional equation $\omega(x;\sigma(z)) = \Gamma_\sigma(z)\omega(x; z)$. Thus we have
\begin{eqnarray*}
&&
Q_n = \sum_{\sigma \in S_n}\, [U\sigma U:U]^{-1} = \frac{w_n(q^{-2})}{(1-q^{-2})^n},\\
&&
\gamma(z) = \prod_{1\leq i < j \leq n}\, \frac{1-q^{z_i-z_j-2}}{1-q^{z_i-z_j}} = \prod_{i < j}\, \frac{q^{z_j}-q^{z_i-2}}{q^{z_j}-q^{z_i}},\\
&&
\Gamma_\sigma(z) = \frac{G_n(\sigma(z))}{G_n(z)}, \quad \mbox{(by Theorem~\ref{th: feq})}, 
\end{eqnarray*}
and \eqref{4-first step} becomes
\begin{eqnarray} \label{second step}
\omega(x;z) &=& 
\frac{(1-q^{-2})^n}{w_n(q^{-2})\cdot G_n(z)} \sum_{\sigma \in S_n}\, 
\sigma\left(\prod_{1\leq i<j\leq n}\, \frac{(q^{z_j}-q^{z_i+1})(q^{z_j}-q^{z_i-2})}{q^{z_j}-q^{z_i}} \delta(x;z)\right),
\end{eqnarray}
where
\begin{eqnarray} \label{delta}
\delta(x; z) = \delta(x; s) = \int_{U} \abs{\bfd(\nu \cdot x)}^s d\nu.
\end{eqnarray}
Hence the problem is reduced to the calculation of $\delta(x;z)$. Let $j = j_n \in K$ be the matrix whose all the anti-diagonal entries are $1$ and other entries are $0$, and set $\check{\pi}^\alp = j \cdot \pi^\alp \in K\cdot \pi^\alp$ for each $\alp \in \Lam_n$, and $jz = (z_n, z_{n-1},\ldots, z_1)$ for $z \in \C^n$. We will prove the next proposition in  \S 4.2. 

\begin{prop} \label{prop: cal of delta}
For any $\alp \in \Lam_n$, one has
\begin{eqnarray*}
\delta(\check{\pi}^\alp; z) =  
\dfrac{c_{odd}(\alp)\cdot q^{\pair{\lam_\alp}{z_0}+\pair{\lam_\alp}{jz}}}{\prod_{\ell \in I_{odd}(\alp)} (q^{z_{n-\ell+1}}-q^{z_{n-\ell}+1})},
\end{eqnarray*}
where $z_0$ is defined in Theorem~\ref{th: explicit}.
\end{prop}

\medskip
Admitting Proposition~\ref{prop: cal of delta} for a while, we substitute the value $\delta(\check{\pi}^\alp; z)$ into \eqref{second step} with $x = \pi^\alp$. Then, for any $\alp \in \Lam_n$, we have
\begin{eqnarray*}
&&
\Psi(\pi^\alp; z) = \omega(\pi^\alp; z) \cdot G_n(z) = \omega(\check{\pi}^\alp;z) \cdot G_n(z)\\
&=&
\frac{(1-q^{-2})^n \cdot c_{odd}(\alp) \cdot q^{\pair{\lam_\alp}{z_0}}}{w_n(q^{-2})} 
\sum_{\sigma \in S_n}\,\sigma\left(  
\frac{q^{\pair{\lam_\alp}{jz}}}{ \prod_{\ell \in I_{odd(\alp)}} (q^{z_{n-\ell+1}}-q^{z_{n-\ell}+1})} 
\prod_{i<j}\frac{(q^{z_j}-q^{z_i+1})(q^{z_j}-q^{z_i-2})}{q^{z_j}-q^{z_i}} \right)\\
&=&
\frac{(1-q^{-2})^n \cdot c_{odd}(\alp) \cdot q^{\pair{\lam_\alp}{z_0}}}{w_n(q^{-2})} 
\sum_{\sigma \in S_n}\,\sigma\left(  
\frac{q^{\pair{\lam_\alp}{z}}}{ \prod_{\ell \in I_{odd(\alp)}} (q^{z_{\ell}}-q^{z_{\ell+1}+1})} 
\prod_{i<j}\frac{(q^{z_i}-q^{z_j+1})(q^{z_i}-q^{z_j-2})}{q^{z_i}-q^{z_j}} \right),
 \end{eqnarray*}
which completes the proof of Theorem~\ref{th: explicit}. \qed
 
\begin{rem} \label{rem: comparison}
{\rm 
When $n=2$, the formula in Theorem~\ref{th: explicit} coincides with that in Proposition~\ref{prop: size2}. For general $n$, we take the main term of $\omega(\pi^\alp; z)$ for each $\alp \in \Lam_n$, and set
\begin{eqnarray} \label{main of sph}
Q(\alp; z) = 
\sum_{\sigma \in S_n}\, 
\sigma \left( \frac{ q^{\pair{\lam_\alp}{z}} }{ \displaystyle{\prod_{\ell \in I_{odd}(\alp)}} (q^{z_{\ell}} - q^{z_{\ell+1}+1 })} \prod_{i<j} \frac{(q^{z_i}-q^{z_j+1})(q^{z_i}-q^{z_j-2}) }{ q^{z_i}-q^{z_j}} \right).
\end{eqnarray}
Then we see $Q(\alp; z)$ is holomorphic for $z \in \C^n$ and linearly independent with respect to $\alp \in \Lam_n$ (cf. Theorem~\ref{th: feq}, Theorem~\ref{th: spTrans}). Corresponding main term of the spherical function on $GL_n(k)$ is a specialization of Hall-Littlewood polynomial 
\begin{eqnarray*}
Q(\lam; z) = \sum_{\sigma \in S_n} \sigma\left(q^{\pair{\lam}{z}} \prod_{i<j} \frac{q^{z_i}-q^{z_j-1} }{ q^{z_i}-q^{z_j}} \right), \quad \lam \in \wt{\Lam_n}.
\end{eqnarray*}
We note here that specializations of Hall-littlewood polynomials also appear in the main term of spherical functions on alternating forms ($X_n \subset GL_{2n}(k)$) as 
\begin{eqnarray*}
Q^{(A)}(\lam; z) = \sum_{\sigma \in S_n} \sigma\left(q^{\pair{\lam}{z}} \prod_{i<j} \frac{q^{z_i}-q^{z_j-2} }{ q^{z_i}-q^{z_j}} \right), \quad \lam \in \wt{\Lam}_n,
\end{eqnarray*}
and on unramified hermitian forms ($X_n \subset GL_n(k'), \; k'/k$ is unramified quadratic ) as
\begin{eqnarray*}
Q^{(H)}(\lam; z) = \sum_{\sigma \in S_n} \sigma\left(q^{\pair{\lam}{z}} \prod_{i<j} \frac{q^{z_i}+q^{z_j-1} }{ q^{z_i}-q^{z_j}} \right), \quad \lam \in \wt{\Lam}_n,
\end{eqnarray*}
(cf. \cite{JMSJ}, \cite{French}, \cite{HS1}, \cite{Mac}). 
In the present case, the shape of $Q(\alp; z)$ is quite different from them. For $n=2$ and even $\alp$, $Q(\alp; z)$ has a relation to Askey-Wilson polynomials, for mode details, see Remark~\ref{rem: A-W}.
} 
\end{rem}

\slit
{\bf 4.2.}\; In this subsection we prove Proposition~\ref{prop: cal of delta}. 
We decompose $U = (U \cap B)U_1$ with 
\begin{eqnarray*}
U_1 = \set{\nu \in GL_n(\calO)}{\begin{array}{ll}
\nu_{ii} = 1 & \mbox{for } 1 \leq i \leq n\\
\nu_{ij} = 0 & \mbox{for } 1 \leq j < i \leq n\\
\nu_{ij} \in \calP & \mbox{for } 1 \leq i < j \leq n\\
\end{array}}.
\end{eqnarray*}
Since $\abs{\bfd(\nu\cdot x)}^s = \abs{\bfd(x)}^s$ for $\nu \in U \cap B$, we see $\delta(x;s) = \int_{U_1}\abs{\bfd(\nu\cdot x)}^s d\nu$.
We have only to consider $\check{\pi}^\alp$ for $\alp \in \Lam_n^+$, since 
$$
\delta(\check{\pi}^{\alp+(2e)}; z) = q^{-e(s_1+2s_2+\cdots + ns_n)}\delta(\check{\pi}^\alp; s) = q^{e(z_1+\cdots + z_n)}\delta(\check{\pi}^\alp; z), \quad e \in \Z.
$$

\begin{lem} \label{lem: entries}
Let $\alp = (\alp_i) \in \Lam_n^+$ and $1 \leq i \leq n$, and assume $(\alp_{n-i+2}, \alp_{n-i+3},\ldots, \alp_n) \in \Lam_{i-1}^+$ if $i >1$. For $\nu \in U_1$, denote by $c(i,j)$ the $(i,j)$-entry of $\nu \cdot \check{\pi}^\alp$.

\mmslit{\rm (1)} If $\alp_{n-i+1}$  is even, say $2e$, then
$c(i, i) = \pi^e$ and $c(i,j) \in \calP^{2e+1}$ for $1 \leq j < i$.

\mmslit{\rm (2)} Assume $\alp_{n-i+1} = \alp_{n-i}$ are odd, say $2e-1$. Then
\begin{eqnarray*}
&&
\twomatrix{c(i,i)}{c(i,i+1)}{c(i+1,i)}{c(i+1,i+1)} = \twomatrix{\pi^e\Trd(u)}{-\varPi^{2e-1}}{\varPi^{2e-1}}{0},\\
&&
c(i,j), c(i+1,j) \in \calP^{2e} \mbox{ for } 1 \leq j \leq i-1,
\end{eqnarray*}
where $u\varPi$ is the $(i,i+1)$-entry of $\nu$.
\end{lem}

\proof 
We see the results by a direct calculation of $\nu\check{\pi}^\alp \times \nu^*$, where we notice that diagonal entries belong to $k$. 
\qed

\begin{lem}  \label{lem: decompose}
Let $\alp = (\alp_i) \in \Lam_n^+$ and $m \leq n$. Assume that $\beta = (\alp_{n-m+1}, \alp_{n-m+2},\ldots,\alp_n) \in \Lam_m^+$. Then for any $\nu \in U_1$, it holds 
\begin{eqnarray*}
\Nrd((\nu\cdot \check{\pi}^\alp)^{(m)}) = \Nrd(\check{\pi}^\beta), \quad
v_\pi(d_m(\nu\cdot \check{\pi}^\alp)) = v_\pi(d_m(\check{\pi}^\alp)) = \frac12\abs{\beta}
\end{eqnarray*}
where $v_\pi(\;)$ is the additive value on $k$.
\end{lem}

\proof
By using Lemma~\ref{lem: entries} consecutively for $i \geq 1$, we see for $\nu \in U_1$, 
\begin{eqnarray*}
&&
d_m(\nu \cdot \check{\pi}^\alp)^2 = \Nrd((\check{\pi}^\alp)^{(m)}) = \Nrd(\check{\pi}^\beta)\; \mbox{and}\;
v_\pi(d_m(\nu\cdot \check{\pi}^\alp)) = \frac12 \abs{\beta}. \hspace{3.5cm} \qed
\end{eqnarray*}

\begin{prop}   \label{prop: even}
If $\alp$ is even, then $\lam_\alp = \frac12\alp$ and one has
\begin{eqnarray*}
\delta(\check{\pi}^\alp;z) = q^{\pair{\lam_\alp}{jz}+\pair{\lam_\alp}{z_0}},
\end{eqnarray*}
where $z_0$ is defined in Theorem~\ref{th: explicit}.
\end{prop}

\proof
Assume that $\alp$ is even and write $\lam = \frac12 \alp = (\lam_i)$. Since we may apply Lemma~\ref{lem: decompose} for every $m$, we have
\begin{eqnarray*}
\delta(\check{\pi}^\alp;z) &= & \prod_{i=1}^n \abs{d_i(\check{\pi}^\alp)}^{s_i}
=
\prod_{i=1}^n q^{-(\lam_n+\lam_{n-1}+\cdots +\lam_{n-i+1})s_i}\\
&=&
\prod_{i=1}^n q^{-\lam_i(s_{n-i+1}+\cdots +s_{n-1}+s_n)}
=
q^{\sum_i(\lam_i(z_{n-i+1}-n+2i-1}\\
&=&
q^{\pair{\lam}{jz}+\pair{\lam}{z_0}}. \hspace{8cm}\qed
\end{eqnarray*}


\medskip
\begin{lem} \label{lem: odd step1}
Let $\alp = (\gamma, \beta) \in \Lam_n^+$ and $\beta \in \Lam_m^+$ and assume $\beta_1=\beta_2=2e-1$. For each $\nu \in U_1$, denote by $c_\nu$ the $(m,m)$-entry of the inverse of the $m\times m$-block of $\nu\cdot \check{\pi}^\alp$. Then one has
\begin{eqnarray}
&&  \label{lem: v_pi eq}
v_\pi(d_{m-1}(\nu\cdot \check{\pi}^\alp)) = \frac12\abs{\beta}+v_\pi(c_\nu),\\
&& \label{volume of c_nu}
vol(\set{\nu \in U_1}{v_\pi(c_\nu)=r}) = (1-q^{-1})q^{-r-e+1}, \quad \mbox{for } r \geq -e+1,
\end{eqnarray}
which depends only on the choice of first $m$ rows of $\nu$. Here we normalize the measure on $U_1$ as $vol(U_1) = 1$.
\end{lem}

%
\proof
By \eqref{miner det}, Lemma~\ref{lem: decompose}, and the fact $c_\nu \in k$,   
we have
\begin{eqnarray*}
d_{m-1}(\nu\cdot \check{\pi}^\alp)^2 
&=& \Nrd((\nu\cdot \check{\pi}^\alp)^{(m-1)}) = \Nrd((\nu\cdot \check{\pi}^\alp)^{(m)})\Nrd(c_\nu) = \Nrd(\check{\pi}^\beta)\Nrd(c_\nu)\\
&=& \Nrd(\check{\pi}^\beta)\, c_\nu^2.
\end{eqnarray*}
Then we obtain the identity \eqref{lem: v_pi eq} by \eqref{Nrd-alp}.
We decompose $\nu \in U_1=U_1(n)$ as 
\begin{eqnarray*}
\nu = \twomatrix{\nu_1}{\nu_1w}{0}{\nu_2}, \quad \nu_1 \in U_1(m), \; \nu_2 \in U_1(n-m), \; w \in M_{m, n-m}(\calP), 
\end{eqnarray*}
then
\begin{eqnarray}  \label{miner m}
(\nu\cdot \check{\pi}^\alp)^{(m)} = \nu_1(\check{\pi}^\beta+w\cdot \check{\pi}^\gamma)\nu_1^*.
\end{eqnarray}
For the simplicity of notation, we set $\beta = (b_m, b_{m-1}, \ldots, b_1)$, where  $2e-1=b_m=b_{m-1} (= b, \; \mbox{say})$ by the assumption.
Set $M = \check{\pi}^\beta+w\cdot \check{\pi}^\gamma \in X_m$. Then we may decompose $M$ as follows:
\begin{eqnarray}
&&
M = M_1 \times M_2, \nonumber \\
&& \label{decomp of M}
M_1= Diag(\varPi^{b_1}u_1,\ldots,\varPi^{b_{m-2}}u_{m-2}, -\varPi^{b}u^*, u\varPi^{b}), \quad u_i, u \in \calO^\times,
\end{eqnarray}
and $M_2 \in K_m$ is congruent modulo $\calP^2$ to the matrix which has diagonal entry $1$ associated to even $b_i$ and diagonal block $H_1$ associated to an odd pair $b_i, b_{i+1}$ and all the other entries are $0$, where $H_1=\twomatrix{0}{1}{1}{0}$.
For example, if $\beta = (3,3,1,1,0,0,0) \in \Lam_7$ then $M_2 \equiv \begin{pmatrix}1_3&0&0\\0&H_1&0\\0&0&H_1\end{pmatrix} \pmod{\calP^2}$.
Since any entry of $w\cdot \check{\pi}^\gamma$ belongs to $\calP^{b+2}$, we see 
\begin{eqnarray}  \label{lowest 2}
(\mbox{the lowest two rows of } M^{-1}) \equiv 
\begin{pmatrix}
0& \cdots & 0 & 0& \varPi^{-b}u^{-1}\\
0& \cdots & 0 & -u^{* -1}\varPi^{-b} & 0
\end{pmatrix}\pmod{\calP^{-b+2}}.
\end{eqnarray}
Denote by $\varPi v^*$ the $(m-1,m)$-entry of $\nu_1 \in U_1(m)$, where $v \in \calO$. Then the $(m,m-1)$-entry of $\nu_1^{* -1}$ is $v\varPi$, and we see by \eqref{lowest 2}
\begin{eqnarray*}
c_\nu &=& (x_1, x_2, \ldots, x_{m-2}, -u^{*-1}\varPi^{-b}+x_{m-1}, v\varPi^{-b+1}u^{-1}+x_m) \begin{pmatrix}y_1\\\vdots\\y_{m-2}\\-\varPi v^*\\1\end{pmatrix}\\
 &\equiv& 
-u^{* -1}\varPi^{-b+1}v^* +v\varPi^{-b+1}u^{-1} \equiv \pi^{-e}\Trd(vu^{-1}) \pmod{\calP^{-b+2}},
\end{eqnarray*}
where $x_1, \ldots, x_{m-2} \in \calP^{-b+1}$, $x_{m-1}, x_m \in \calP^{-b+3}$, $y_1, \ldots, y_{m-2}\in \calP$, they are determined by $\nu_1$ and $w$, and the above column vector is the $m$-th column of $\nu_1$. Since $c_\nu \in k$, we may write  
\begin{eqnarray*}
c_\nu = \pi^{-e+1}(\Trd(vu^{-1}) + \pi z), \quad z \in \calo.
\end{eqnarray*}
Here $z$ is determined by $\alp$ and $w$ and $\nu_1$ except the $(m-1,m)$-entry $\varPi v^*$, and $u \in \calO^\times$ is determined by $\alp$ and $w$ (cf. \eqref{decomp of M}). Hence $v_\pi(c_\nu) \geq -e+1$, which is determined by the choice of $v \in \calO$ and independent of the choice of $\nu_2$. For $r \geq -e+1$,
\begin{eqnarray*}
vol(\set{\nu \in U_1}{v_\pi(c_\nu) = r}) &=& vol\set{v \in \calO}{\Trd(v) \in \pi^{r-e+1}\calo^\times}\\
&=&
q^{-4\ell}\sharp\set{\ol{u} \in \calO/\calP^{2\ell}}{\Trd(u) \in \pi^{r+e-1}\calo^\times}\\
& = &
q^{-r-e+1}(1-q^{-1}),
\end{eqnarray*}   
where we take $\ell \geq r+e$ and we regard $\Trd$ as $(q^{3\ell}:1)$-mapping from $\calO/\calP^{2\ell}$ onto $\calo/\frp^\ell$.
\qed

\bigskip
Now, we assume $\alp$ is has odd entries and set $I_{odd}(\alp) = \{\ell_1, \ldots, \ell_k\}$ (cf. \eqref{I_odd}) and $\alp_{\ell_j} = 2e_j-1$. 
For $m=n-\ell_j,\; 1 \leq j \leq k$, we use Lemma~\ref{lem: odd step1}, and for the other $m$, we may use Lemma~\ref{lem: decompose}. Then we obtain
\begin{eqnarray*}
\delta(\check{\pi}^\alp;z) &=& 
\prod_{m \ne n-\ell_j} q^{-\frac12(\alp_{n-m+1}+\cdots +\alp_n)s_m}\\
&&\times
\prod_{j=1}^k  \sum_{r_j \geq -e_j+1}\, (1-q^{-1}) q^{-(r_j+e_j-1)}q^{-(\frac12(\alp_{\ell_j}+\alp_{\ell_j+1}+\cdots +\alp_n)+r_j)s_{n-\ell_j}}\\
&=&
\prod_{m \ne n-\ell_j} q^{-\frac12(\alp_{n-m+1}+\cdots +\alp_n)s_m}\times
\prod_{j=1}^k \frac{1-q^{-1}}{1-q^{-1-s_{n-\ell_j}}}\cdot q^{-\frac12(1+\alp_{\ell_j+1}+\alp_{\ell_j+2}+\cdots+\alp_n)s_{n-\ell_j}}\\
&&\hspace*{10cm}
(\alp_{\ell_j}=\alp_{\ell_j+1}=2e_j-1)\\
&=&
\prod_{i=1}^n q^{-\frac12 \alp_i(s_{n-i+1}+\cdots +s_n)} \times
\prod_{j=1}^k \frac{(1-q^{-1})q^{\frac12 (z_{n-\ell_j}-z_{n-\ell_j+1}+2) }}{1-q^{z_{n-\ell_j}-z_{n-\ell_j+1}+1}}\\
&=&
\prod_{i=1}^n q^{\frac12 \alp_i(z_{n-i+1}-n+2i-1)} \times
\prod_{j=1}^k \frac{ (1-q^{-1})q^{\frac12 (z_{n-\ell_j}+z_{n-\ell_j+1})+1} }
{ q^{z_{n-\ell_j+1}} - q^{z_{n-\ell_j}+1}}\\
&=&
q^{\pair{\frac{\alp}{2}}{z_0}} \times \left\{\prod_{i=1}^n q^{\frac12\alp_i z_{n-i+1}}
\prod_{j=1}^k q^{\frac12 (z_{n-\ell_j}+z_{n-\ell_j+1})}\right\} \times
\prod_{j=1}^k \frac{ (1-q^{-1})q} 
{ (q^{z_{n-\ell_j+1}} - q^{z_{n-\ell_j}+1})}\\
&=&
q^{\pair{\lam_\alp}{z_0}
-\frac12\sum_{j=1}^k(-n+2\ell_j-1)+(-n+2(\ell_j+1)-1) }\cdot q^{\pair{\lam_\alp}{jz}}
\times \prod_{j=1}^k \frac{ (1-q^{-1})q} 
{ (q^{z_{n-\ell_j+1}} - q^{z_{n-\ell_j}+1})}\\
&=&
\frac{(1-q^{-1})^k q^{ \sum_{j=1}^k (n-2\ell_j+1)}}{\prod_{\ell \in I_{odd}(\alp)} (q^{z_{n-\ell+1}}-q^{z_{n-\ell}+1}) } \times
{q^{\pair{\lam_\alp}{z_0}+\pair{\lam_\alp}{jz}}}\\
&=&
\frac{c_{odd}(\alp)}{\prod_{\ell \in I_{odd}(\alp)} (q^{z_{n-\ell+1}}-q^{z_{n-\ell}+1}) } \times
{q^{\pair{\lam_\alp}{z_0}+\pair{\lam_\alp}{jz}}}.
\end{eqnarray*}
\qed

\vspace{1cm}
\Section{Schwartz space $\SKX$}

\mslit
{\bf 5.1.} 
In this subsection, we study $\hec$-module structure of $\SKX$ through the spherical transform $F$, so we need to recall Theorem~\ref{th: spTrans}. 
For each $\alp \in \Lam_n$, we denote by $\vphi_\alp \in \SKX$ the characteristic function of $K\cdot \pi^\alp$ and $\check{\vphi}(x)=\vphi(x^{-1})$ for $\vphi \in \SKX$. Then 
\begin{eqnarray*}
F(\check{\vphi_\alp}) = vol(K\cdot \pi^\alp) \Psi(\pi^\alp;z) \equiv \Psi(\pi^\alp;z),
\end{eqnarray*}
where $\equiv$ means $\equiv({\rm mod\,}\R^\times)$ in this subsection. 
The image $F(\SKX) = \gen{\Psi(\pi^\alp; z): \alp \in \Lam_n}_\C$ is an ideal of $\calR = \C[q^{\pm z_1},\ldots, q^{\pm z_n}]^{S_n}$, where $\calR$ is isomorphic to $\hec$ by Satake isomorphism (cf. \eqref{lam_z}).

\medskip
To make sure of it, we note the results for $n = 1, 2$, which is easily seen by definition of spherical function for $n = 1$ and Proposition~\ref{prop: size2} for $n=2$.

\begin{prop}
Assume $n = 1$ or $2$. Then  the spherical transform $F: \SKX \longrightarrow \calR$ is an $\hec$-module isomorphism, and $\SKX$ is generated as an $\hec$-module by $\vphi_{\bf0}$ for $n=1$ and $\vphi_{(-1,-1)}$ for $n =2$.
\end{prop}

As a corollary we see the following, where we omit the proof since it can be proved similarly and more easily to Proposition~\ref{prop: param of sph 3,4}.   

\begin{prop}
Assume $n = 1$ or $2$. Then any spherical function on $X$ is a constant multiple of $\Psi(x; z_0)$ for some $z_0 \in \C^n/S_n$.
\end{prop} 

\medskip
In the rest of this subsection, we consider the case $n \geq 3$. For simplicity of notation, we set $x_i = q^{z_i}, \; 1 \leq i \leq n$ and denote by $s_i(n)$ the fundamental symmetric polynomial in $x_1, \ldots, x_n$, for $1 \leq i \leq n$. Then
\begin{eqnarray}
\calR = \C[x_1^{\pm1}, \ldots, x_n^{\pm1}]^{S_n} = \C[s_1(n), \ldots, s_n(n), s_n(n)^{-1}] = \calR_0[s_n(n)^{-1}],
\end{eqnarray}
where $\calR_0 = \C[s_1(n), \ldots, s_n(n)]$. We denote by $J$ the ideal of $\calR_0$ generated by the subset $\set{\Psi(\pi^\alp; z)}{\alp \in \Lam_n, \; \alp_n = 0, -1}$ of $\calR_0$. Then, since $\Psi(\pi^{\alp+(2e)};z)=s_n(n)^e\Psi(\pi^\alp; z)$ for $e \in \Z$,  we see that $F(\SKX) = J \tensor_{\calR_0}\calR$.

\mslit
For a (fixed) rational function $c(x)$ of $x_1, \ldots, x_n$ and $\mu \in \Z^n$, we define
\begin{eqnarray}
P(c(x), \mu;x) = \sum_{\sigma \in S_n}\, \sigma(x_1^{\mu_1}\cdots x_n^{\mu_n}\cdot c(x)) \in \C(x_1, \ldots, x_n)^{S_n},
\end{eqnarray}
where $S_n$ acts on $x_1, \ldots, x_n$ by permutation of indices.

\begin{lem} \label{lem: poly-reduction}
For any $\lam \in \wt{\Lam}_n$ with $\lam_n \geq 0$, the rational function $P(c(x), \lam;x)$ belongs to the $\calR_0$-module generated by the set $\set{P(c(x), \mu; x)}{\mu \in \Sigma_n}$, where $\Sigma_n =\\ \set{\mu = (\mu_i) \in \Z^n}{0 \leq \mu_i \leq n-1}$.
\end{lem}

\proof
Since $s_i(n+1) = s_i(n)+s_{i-1}(n)x_{n+1}$, we have the following by induction on $n$.  
\begin{eqnarray} \label{down}
x_\ell^r = \sum_{i=1}^\ell\, (-1)^{i-1}s_i(n)x_\ell^{r-i}, \quad (1 \leq \ell \leq n), \quad \mbox{if }\; r \geq n.
\end{eqnarray}
For each $\lam \in \wt{\Lam}_n$, we set $\lam^{(\ell,i)} = (\lam_1, \ldots,\lam_{\ell-1}, \lam_\ell-i, \lam_{\ell+1},\ldots,\lam_n) \in \Z^n$. Then, by \eqref{down}, we have
\begin{eqnarray}
P(c(x),\lam;x) = \sum_{i=1}^n\, (-1)^i s_i(n)P(c(x),\lam^{(\ell,i)};x), \quad \mbox{if }\; \lam_\ell \geq n.
\end{eqnarray}
Taking this procedure for every $\ell$ with $\lam_\ell \geq n$, we have the result. \qed

\medskip
By computer calculation, it is possible to express symmetric polynomials in terms of $s_i(n), \; 1 \leq i \leq n$, if the variable $n$ and degrees of polynomials are small. The author owes Satoshi Murai for a program using Macaulay2(\cite{Macaulay}), which worked well for size $n \leq 4$. Further it is possible to check for a polynomial whether it is contained in a fixed ideal of $\calR_0$ or not, by Macaulay2.
Thus we have the following proposition.   

\begin{prop} \label{prop: size3}
Assume $n=3$.\\
{\rm (1)} $F(\SKX)$ is an ideal of $\calR$ generated by $\Psi(1_3; z) \equiv s_1s_2-q^2(q^{-2}+q^{-1}+1)^2 s_3$ and $\Psi(\pi^{(0,-1,-1)}; z) \equiv s_1^2- q^2(q^{-2}+q^{-1}+1)^2 s_2$, where $s_i = s_i(3), \; 1 \leq i \leq 3$, and it is non-principal.

\mmslit
{\rm (2)} The $\hec$-module $\SKX$ is generated by $\vphi_{\bf0}$ and $\vphi_{(1,1,0)}$, and it is not monomial.
\end{prop}

\proof
(1) In consideration of the explicit formula of $\Psi(\pi^\alp; z)$ (Theorem~\ref{th: explicit}), we set 
\begin{eqnarray}
&&
c_3(x) = \prod_{1 \leq i < j \leq 3}\frac{(x_i-qx_j)(x_i-q^{-2}x_j)}{x_i-x_j}, \nonumber \\
&&
P1(a,b;x) = P(c_3(x), (a,b,0); x) \equiv \Psi(\pi^{(2a,2b,0)};z), \nonumber \\
&&
P2(a;x) = P(\frac{c_3(x)}{x_2-qx_3}, (a,0,0); x) \equiv \Psi(\pi^{(2a,-1,-1)}; z), \nonumber\\ 
&&
P3(a;x) = P(\frac{c_3(x)}{x_1-qx_2}, (a,a,0); x) \equiv \Psi(\pi^{(2a-1,2a-1,0)};z). 
\end{eqnarray}
Then $J$ is generated by the set $\set{P1(a,b;x)}{a \geq b \geq 0}\cup \set{P2(a;x)}{a \geq 0} \cup  \set{P3(a;x)}{a \geq 1}$.
Let $J_0$ be the ideal of $\calR_0$ generated $P1(0,0;x) \equiv \Psi(\pi^{(0,0,0)}; z) = \Psi(1_3; z)$ and $P2(0;x) \equiv \Psi(\pi^{(0,-1,-1)}; z)$, where their values are calculated as above.  It is clear that $J_0$ is not principal. 
We may check by computer that the set $\set{P1(a,b;x)}{a, b = 0, 1, 2} \cup \set{P2(a;x)}{a = 1, 2} \cup  \set{P3(a;x)}{a = 1, 2 }$ is contained in $J_0$, which brings $J = J_0$ by Lemma~\ref{lem: poly-reduction}.
As for $P3$, the above set is enough, since we may use $x_1^3x_2^3= s_2x_1^2x_2^2 - s_1s_3x_1x_2+s_3^2$ instead of \eqref{down}.
\\
(2) Since $F$ is an injective $\hec$-module map and $F(\check{\vphi_\alp}) \equiv \Psi(\pi^\alp; z)$, we see that $\SKX$ is generated by $\check{\vphi}_{\bf0} = \vphi_{\bf0}$ and $\check{\vphi}_{(0,-1,-1)} = \vphi_{(1,1,0)}$ and is not monomial. 
\qed

\begin{rem}
{\rm In the above, since any polynomial in $\set{P1(a,b; x)}{0 \leq a < b \leq 2}$ is not combined with $\Psi(\pi^\alp; z), \alp \in \Lam_3$, it is not assured to be contained even in $J$. Fortunately, it is contained in $J_0$. }
\end{rem}

In the similar way to size 3, we obtain the following result for size $4$. We have to consider 5 types of polynomials associated with $\Psi(\pi^\alp; z)$ according to the location of odd entries of $\alp \in \Lam_4$(cf. Theorem~\ref{th: explicit}). Then, by Lemma~\ref{lem: poly-reduction}, it is enough to verify a certain finite set of polynomials. In this case, as with the case of size 3,  some polynomials are not combined with $\Psi(\pi^\alp; z)$, but they are contained in the ideal fortunately. The second claim follows from the first one and the fact the spherical transform $F$ is an injective $\hec$-module map.

\begin{prop} \label{prop: size4}
Assume $n=4$.\\  
{\rm (1)} $F(\SKX)$ is an ideal of $\calR$ generated by two elements
\begin{eqnarray*}
&&
\Psi(1_4; z) 
\equiv
s_1 s_2 s_3 - q^2(q^{-2}+q^{-1}+1)^2 s_1^2 s_4 \\
&&
\hspace*{3cm}  -q^2(q^{-2}+q^{-1}+1)^2 s_3^2 +q^3(q^{-2}+1)(q^{-1}+1)^4 s_2s_4, \\
&&
\Psi(\pi^{(-1,-1,-1,-1)}; z)\equiv s_2^2 -q(q^{-2}+q^{-1}+1)s_1 s_3 + q^3(q^{-2}+1)^2(q^{-2}+q^{-1}+1) s_4,
\end{eqnarray*}
where $s_i = s_i(4), \; 1 \leq i \leq 4$, and it is non-principal.

\mmslit
{\rm (2)} The $\hec$-module $\SKX$ is generated by $\vphi_{\bf0}$ and $\vphi_{(1,1,1,1)}$, and not monimial.\\
 \end{prop}

\medskip
As a corollary of Propositions~\ref{prop: size3} and \ref{prop: size4}, we have the following.

\begin{prop} \label{prop: param of sph 3,4}
Assume $n=3, 4$, and set $\beta = (1,1,0)$ for $n=3$ and  $\beta =(1,1,1,1)$ for $n = 4$.
If $\Psi(1_n; z_0) \ne 0$ and $\Psi(\pi^\beta; z_0) \ne 0$, then any spherical function on $X$ corresponding to $z_0$ is a constant multiple of $\Psi(x; z_0)$.
\end{prop}

\proof
Any spherical function is associated with some $z_0 \in \C^n$ by its eigenvalue (cf. the comment at the end of \S 1). We introduce the pairing on $\SKX \times \CKX$ by 
\begin{eqnarray} \label{pairingSC}
\pair{\vphi}{\Phi} = \int_X \vphi(x)\Phi(x)dx, \quad \vphi \in \SKX, \; \Phi \in \CKX
\end{eqnarray}
where $dx$ is the $G$-invariant measure. Then it satisfies for any $f \in \hec$,
\begin{eqnarray} \label{hec-act}
\pair{f*\vphi}{\Phi} = \pair{\vphi}{\check{f}*\Phi}, \quad \check{f}(g) = f(g^{-1}).
\end{eqnarray}
Assume $\Phi$ is a spherical function on $X$ corresponding to $z_0$ which satisfies the assumption above. 
Denote $\vphi_0 = \vphi_{\bf0}$ and $\vphi_1 = \vphi_\beta$. Then for any $f \in \hec$ and $i = 0, 1$, one has by \eqref{pairingSC} and \eqref{hec-act}
\begin{eqnarray*}
&&
\pair{f*\vphi_i}{\Phi} = \lam_{z_0}(\check{f})\pair{\vphi_i}{\Phi},\\
&&
\pair{f*\vphi_i}{\Psi(\; ; z_0)} = \lam_{z_0}(\check{f})\pair{\vphi_i}{\Psi(\; ; z_0)},
\end{eqnarray*}
where $\pair{\vphi_i}{\Psi(\; ; z_0)} \ne 0$ by the coice of $z_0$. 
Thus 
\begin{eqnarray}   \label{eigenvalue}
\pair{f*\vphi_i}{\Phi} = \frac{\pair{\vphi_i}{\Phi}}{\pair{\vphi_i}{\Psi(\; ;z_0)}} 
\pair{f*\vphi_i}{\Psi(\; ; z_0)}, \quad f \in \hec, \; i=0,1.
\end{eqnarray}
On the other hand, by the commutative diagram \eqref{comm diag}, there is some $g_i \in \hec, \; i = 0,1$ such that
\begin{eqnarray*}
&&
\lam_z(g_1)F(\vphi_0) = \lam_z(g_0)F(\vphi_1) \ne 0 \nonumber \\
\end{eqnarray*}
thus it holds
\begin{eqnarray} \label{equality}
&& 
g_1*\vphi_0 = g_0*\vphi_1 (\ne 0),
\end{eqnarray}
hence
\begin{eqnarray} \label{equality2}
\pair{g_1*\vphi_0}{\Psi(\; ;z_0)} = \pair{g_0*\vphi_1}{\Psi(\; ;  z_0)} \; (\ne 0),
\end{eqnarray}
By \eqref{eigenvalue}, \eqref{equality} and \eqref{equality2}, one sees
\begin{eqnarray*}
c_\Phi:= \frac{\pair{\vphi_0}{\Phi}}{\pair{\vphi_0}{\Psi(\; ;z_0)}} = \frac{\pair{\vphi_1}{\Phi}}{\pair{\vphi_1}{\Psi(\; ;z_0)}} (\ne 0).
\end{eqnarray*}
 Then, since $\SKX = \hec*\vphi_0 + \hec*\vphi_1$, we have
\begin{eqnarray*}
\pair{\vphi}{\Phi} = c_\Phi \pair{\vphi}{\Psi(\; ;z_0)} = \pair{\vphi}{c_\Phi \Psi(\; ;z_0)},  \quad \vphi \in \SKX,
\end{eqnarray*}
which yields $\Phi(x) = c_\Phi\Psi(x; z_0)$ in $\CKX$ as required.
\qed
  
\begin{rem}
{\rm 
It is expected that, for general size $n \geq 5$, the $\hec$-module $\SKX$ is generated by $\vphi_{\bf0}$ and $\vphi_{\beta}$, where $\beta = (1,\ldots,1)$ or $(1,\ldots,1,0)$ according to the parity of $n$, and is not monimial. In other words, it is expected that the ideal $F(\SKX)$ of $\calR$ is generated by $\Psi(\pi^{\bf0};z) = \Psi(1_n; z)$ and $\Psi(\pi^{\beta'};z)$ where $\beta' = (-1,\dots,-1)$ or $(0,-1,\ldots,-1)$, and is not principal.
If this is true for $n$, then the parallel result for $n$ to Proposition~\ref{prop: param of sph 3,4} holds.   
}
\end{rem}

\vspace{2cm}

\noindent
{\bf 5.2.} We introduce the Plancherel formula for size $2$ proved by Y.~Komori. Throughout this subsection we only consider the case of size $2$, hence $X \subset G = GL_2(D)$.  

\mslit
{\bf 5.2.1} Since $G_n$-invariant measure $dx$ on $X_n$ of size $n$ is determined up to constant by the differential form 
\begin{eqnarray}
\dfrac{\wedge_{i=1}^n dx_{ii} \wedge \wedge_{1\leq i<j\leq n}\, dx_{ij}}{\abs{\Nrd(x)}^{\frac{2n-1}{n}}}, \quad (x = (x_{ij})\in X,\; x_{ii} \in k, \; x_{ij} \in D),
\end{eqnarray}
volume $v(K\cdot \pi^\alp)$ in the case of size $2$ is a constant multiple of $q^{\frac{3}{2}\abs{\alp}}/\mu(\pi^\alp, \pi^\alp)$, and by Theorem~\ref{th: vol} we have
\begin{eqnarray}
\mu(\pi^\alp, \pi^\alp) = \left\{\begin{array}{ll}
q^{6\lam_1}(1+q^{-1})(1-q^{-2}) &\mbox{if }\; \alp = (2\lam_1, 2\lam_1)\\
q^{\lam_1+5\lam_2}(1+q^{-1})^2 &\mbox{if }\; \alp = (2\lam_1, 2\lam_2), \; \lam_1 > \lam_2\\
q^{6e-2}(1-q^{-4}) & \mbox{if }\; \alp = (2e-1,2e-1).
\end{array}\right.
\end{eqnarray}
We normalize $dx$ as $v(K\cdot 1_n) = 1$. 
Then, for the characteristic function $\vphi_\alp$ of $K \cdot \pi^\alp, \; \alp \in \Lam_2$ and $\check{\vphi_\alp}(x)=\vphi_\alp(x^{-1})$, we see 
\begin{eqnarray}
&&
\int_X \check{\vphi_\alp}(x)\ol{\check{\vphi_\beta}(x)}dx = \delta_{\alp,\beta}\times v(K\cdot \pi^{(-\alp_2, -\alp_1)}) = \delta_{\alp,\beta}\times v(K\cdot \pi^\alp)  \nonumber\\
&& \label{side X}
= \delta_{\alp, \beta} \times \left\{\begin{array}{ll}
1  &\mbox{if }\; \alp_1 = \alp_2 \in 2\Z\\[1.5mm]
q^{2(\lam_1-\lam_2)}(1-q^{-1})  &\mbox{if }\; \alp = (2\lam_1, 2\lam_2), \; \lam_1 > \lam_2\\
\dfrac{q^{-1}(1+q^{-1})}{1+q^{-2}} &  \mbox{if }\; \alp_1 = \alp_2 \notin 2\Z 
\end{array}\right\}. 
\end{eqnarray}
On the other hand,  
by the definition of $F$ (cf. \eqref{def G_n(z)}, \eqref{sph transf-s}) and Proposition~\ref{prop: size2}, we have
\begin{eqnarray}
F(\check{\vphi_\alp}) &=& v(K\cdot \pi^\alp)\,\Psi(\pi^\alp;z) = v(K\cdot \pi^\alp)\omega(\pi^\alp;z)(q^{z_2}-q^{z_1+1})\nonumber\\[2mm]
&=&  \label{F(alp)}
\left\{\begin{array}{ll}
\frac{1}{1+q^{-2}}\, Q_{\lam}(z)
 & \mbox{if\;} \alp = 2\lam, \lam_1=\lam_2,\\[2mm]
\dfrac{ q^{\lam_1-\lam_2}(1-q^{-1})}{1+q^{-2}}\, Q_\lam(z)
& \mbox{if\;}\alp = 2\lam, \lam_1>\lam_2,\\[2mm]
\dfrac{1-q^{-2}}{1+q^{-2}}\, q^{e(z_1+z_2)} & \mbox{if\;}\alp=(2e-1,2e-1).
\end{array}\right.
\end{eqnarray}
where
\begin{eqnarray}
Q_\lam(z) = (q^{z_1}+q^{z_2})\sum_{\sigma\in S_2} \sigma\left(q^{\lam_1z_1+\lam_2z_2}\frac{(1-q^{z_2-z_1+1})(1-q^{z_2-z_1-2})}{1-q^{2(z_2-z_1)}}\right).
\end{eqnarray}

\mslit
{\bf 5.2.2} Fix $u_i$ as $0< u_i < 1, \; i = 1,2$ and set 
\begin{eqnarray}
&&
H_\ell(y)=\sum_{\sigma \in S_2}\, \sigma\left(q^{-\ell y}\frac{(1-u_1q^{2y})(1-u_2q^{2y})}{1-q^{2y}}\right), \quad (\ell \in \N), \\[2mm]
&&
w(y)=\frac{1-q^{2y}}{(1-u_1q^{2y})(1-u_2q^{2y})}\cdot\frac{1-q^{-2y}}{(1-u_1q^{-2y})(1-u_2q^{-2y})}.
\end{eqnarray}
Take $U = \set{y = \sqrt{-1}t}{0 \leq t \leq 2\pi\log q}$ and a (suitably normalized) measure $dy$ on $U$, one has
\begin{eqnarray}
&& \label{int H_m and H_n}
\int_U H_\ell(y)\ol{H_m(y)}w(y)dy = \left\{\begin{array}{ll}
(1-u_1u_2) & (\ell=m=1),\\
1 &(\ell=m>1),\\
0 & (\ell\ne m),  
\end{array}\right.\\
&& \label{int H_n}
\int_{U} H_\ell(y)w(y)dy=0 \qquad(\ell\geq 1),\\
&& \label{w(y)itself}
\int_{U}w(y)dy = \frac{1}{(1+u_1)(1+u_2)(1-u_1u_2)},
\end{eqnarray}
where $\ol{H_m(y)}$ is the complex conjugate of $H_m(y)$. 
As for \eqref{w(y)itself}, we calculate the integral
\begin{eqnarray*}
\frac{1}{2\pi\sqrt{-1}}\int_{\abs{Y}=1}\frac{1-Y^{-2}}{(1-u_1Y^{-2})(1-u_2Y^{-2})}\cdot\frac{1-Y^2}{(1-u_1Y^2)(1-u_2Y^2)}
\frac{dY}{Y} = \frac{2}{(1+u_1)(1+u_2)(1-u_1u_2)}.
\end{eqnarray*}

\mslit
\begin{rem} \label{rem: A-W}
{\rm 
 The set $\{H_m\}$ essentially coincides with a special case of the Hall-Littlewood limit of the Askey-Wilson polynomials
  \cite{KS96}, that is, the limit $q\to 0$ in the context of $q$-orthogonal polynomials. This follows from the fact that the Hall-Littlewood limit of the Askey-Wilson polynomials
  satisfies the orthogonal conditions (5.23) and (5.24) and that such polynomials with the leading terms $Y^m$ are unique.
}
\end{rem}

\medskip
Set $x = \frac{z_2+z_1}{2}, \; y = \frac{z_2-z_1}{2}$ and $u_1 = q, \; u_2 = q^{-2}$. Then for $\lam = (\lam_1, \lam_2) \in \wt{\Lam}_2$, we have
\begin{eqnarray}
Q_\lam(z) 
&=&
q^{(\lam_1+\lam_2+1)x}(q^y+q^{-y})\sum_{\sigma \in S_2} \sigma \left(q^{-(\lam_1-\lam_2)y}\frac{(1-u_1q^{2y})(1-u_2q^{2y})}{1-q^{4y}}\right)\nonumber \\
&=& 
q^{(\abs{\lam}+1)x}\sum_{\sigma \in S_2} \sigma\left(q^{-(\lam_1-\lam_2+1)y}\frac{(1-u_1q^{2y})(1-u_2q^{2y})}{1-q^{2y}}\right) \nonumber \\
&=& \label{Q_lam}
q^{(\abs{\lam}+1)x}H_{\lam_1-\lam_2+1}(y),
\end{eqnarray}
where we set $\abs{\lam}=\lam_1+\lam_2$.
For $e \in \Z$, we define
\begin{eqnarray}
R_e(x,y) = q^{(e+1)(z_1+z_2)} = q^{2(e+1)x}.
\end{eqnarray}

Now, for a while, we consider $u_1, u_2$ are independent of $q$ and still keep the condition $0<u_1, u_2 < 1$. We set $T = \set{x=\sqrt{-1}t}{0 \leq t \leq 2\pi\log q}$ and $U = \set{y=\sqrt{-1}t}{0 \leq t \leq 2\pi\log q}$,  
and define the inner product on $\C[q^x,q^{-x}, q^y+q^{-y}]$ by
\begin{eqnarray} \label{inner product}
\pair{f}{g} = \int_{T}dx\int_{U}f\ol{g}w(y)dy.
\end{eqnarray}
Then we see, (cf. \eqref{int H_m and H_n} also)
\begin{eqnarray}
&&
\pair{R_e}{R_{e'}} = 0, \quad \mbox{for } e,e' \in \Z, \; e \ne e',\\
&&
\pair{Q_\lam}{Q_\mu} = 0, \quad \mbox{for } \lam, \mu \in \wt{\Lam}_2, \; \lam \ne \mu.  
\end{eqnarray}
If $\abs{\lam} = \lam_1+\lam_2$ is even, by the integral with respect to $x$, we see $\pair{Q_\lam}{H_e} = 0$.
Assume $\abs{\lam}$ is odd. Then $\pair{Q_\lam}{H_e} = 0$ unless $\abs{\lam}= 2e+1$. When $\abs{\lam}=2e+1$, $\lam_1-\lam_2=2n-1 > 0$, and $Q_\lam(x,y) = q^{(\abs{\lam}+1)x}H_{2n}(y)$. Then $\pair{Q_\lam}{R_e} = 0$ by \eqref{int H_n}.

\medskip
We have to consider the case $u_1=q$ and $u_2=q^{-2}$. We fix $u_2=q^{-2}$ and change $u_1$ continuously from $0<u_1<1$ to 
$q$. When $0 < u_1 <1$, the poles of $w(y)$ are, written by variable $Y = q^y$, 
\begin{eqnarray}
&&
Y = \pm\sqrt{u_1}, \; \pm \sqrt{u_2} = \pm q^{-1} \; \mbox{ within $\abs{Y} < 1$},\\
&&
Y = \pm\sqrt{u_1^{-1}}, \; \pm \sqrt{u_2^{-1}} = \pm q \; \mbox{ outside of $\abs{Y} = 1$}.
\end{eqnarray}
According to the change of $u_1$, we change the integration path $\abs{Y} = 1$ as the path does not change the way around these poles (cf. Figures 1,2). 

\begin{figure}[h]
  \begin{minipage}{0.5\linewidth}
    \begin{center}
        \begin{tikzpicture}[scale=0.65]
  
  \draw [-latex] (0,-4)--(0,4) ;
  \draw [-latex] (-6,0)--(6,0) ;

  \filldraw (1,0) circle [radius=0.05] node [above] {$q^{-1}$};
  \filldraw (5,0) circle [radius=0.05] node [above] {$q$};
  \filldraw (-1,0) circle [radius=0.05] node [above] {$-q^{-1}$};
  \filldraw (-5,0) circle [radius=0.05] node [above] {$-q$};

  \filldraw (2,0) circle [radius=0.05] node [below] {$\sqrt{u_1}$};
  \filldraw (4,0) circle [radius=0.05] node [below] {$\sqrt{u_1}^{-1}$};
  \filldraw (-2,0) circle [radius=0.05] node [below] {$-\sqrt{u_1}$};
  \filldraw (-4,0) circle [radius=0.05] node [below] {$-\sqrt{u_1}^{-1}$};

  \draw [ultra thick, decoration={markings, 
    mark=at position 0.2 with {\arrow{latex}}
  },
  postaction={decorate}
  ] (0,0) circle [radius=3];

\end{tikzpicture}
\end{center}
\caption{the case $0< u_1 <1$}
  \label{fig:a1}
\end{minipage}
\begin{minipage}{0.5\linewidth}
  \begin{center}
    \begin{tikzpicture}[scale=0.65]
  
  \draw [-latex] (0,-4)--(0,4) ;
  \draw [-latex] (-6,0)--(6,0) ;

  \filldraw (1,0) circle [radius=0.05] node [above] {$q^{-1}$};
  \filldraw (5,0) circle [radius=0.05] node [above] {$q$};
  \filldraw (-1,0) circle [radius=0.05] node [above] {$-q^{-1}$};
  \filldraw (-5,0) circle [radius=0.05] node [above] {$-q$};

  \filldraw (2,0) circle [radius=0.05] node [below=6] {$\sqrt{q}^{-1}$};
  \filldraw (4,0) circle [radius=0.05] node [below=6] {$\sqrt{q}$};
  \filldraw (-2,0) circle [radius=0.05] node [below=6] {$-\sqrt{q}^{-1}$};
  \filldraw (-4,0) circle [radius=0.05] node [below=6] {$-\sqrt{q}$};

  \draw [ultra thick,decoration={markings, 
    mark=at position 0.2 with {\arrow{latex}}
  },
  postaction={decorate}
  ] (10:3) arc [start angle=10, end angle=170,radius=3]
  ;
  \coordinate (a) at ($(-4,0)!1!30:(-3.5,0)$) ;
  \coordinate (b) at ($(-2,0)!1!-60:(-2.5,0)$) ;
  \coordinate (c) at ($(4,0)!1!30:(3.5,0)$) ;
  \coordinate (d) at ($(2,0)!1!-60:(2.5,0)$) ;

  \draw [ultra thick] (170:3) .. controls (-3,0) and (-3.4,0) .. (a) arc [start angle=30, end angle=300,radius=0.5] -- (b) arc [start angle=120, end angle=-150,radius=0.5] .. controls (-2.6,0) and (-3,0) .. (190:3) ;
  \draw [ultra thick] (190:3) arc [start angle=190, end angle=350,radius=3] ;

  \draw [ultra thick] (-10:3) .. controls (3,0) and (3.4,0) .. (c) arc [end angle=120, start angle=-150,radius=0.5] -- (d)  arc [end angle=30, start angle=300,radius=0.5]  .. controls (2.6,0) and (3,0) .. (10:3) ;

\end{tikzpicture}
\end{center}
\caption{the case $u_1=q>1$}
\end{minipage}
\end{figure}

\newpage
\noindent
{\bf 5.2.3} We calculate the norm $\pair{F(\check{\vphi_\alp})}{F(\check{\vphi_\alp})}$ by using \eqref{F(alp)}, \eqref{int H_m and H_n}, \eqref{w(y)itself}, and \eqref{Q_lam}.
When $\alp = (2e,2e)$, 
\begin{eqnarray}
\pair{F(\check{\vphi_\alp})}{F(\check{\vphi_\alp})}= \frac{(1-q^{-1})}{(1+q^{-2})^2}, \qquad(v(K\cdot \pi^\alp) = 1);
 \end{eqnarray}
when $\alp = 2\lam, \; \lam_1>\lam_2$, the value is 
 \begin{eqnarray}
\frac{q^{2(\lam_1-\lam_2)}(1-q^{-1})^2}{(1+q^{-2})^2} = v(K\cdot \pi^\alp)\times \frac{1-q^{-1}}{(1+q^{-2})^2};
\end{eqnarray}
when $\alp = (2e-1, 2e-1)$, the value is
 \begin{eqnarray}
\frac{(1-q^{-2})^2}{(1+q^{-2})^2 (1+q)(1+q^{-2})(1-q^{-1}) } = \frac{q^{-1}(1-q^{-2})}{(1+q^{-2})^3} = v(K\cdot \pi^\alp) \times \frac{1-q^{-1}}{(1+q^{-2})^2} .
     \end{eqnarray}

\mslit
By comparison with \eqref{side X}, we normalize the inner product \eqref{inner product} by multiplying $(1+q^{-2})^2/(1-q^{-1})$ and keep the notation, then we obtain 
\begin{eqnarray}
\int_X \check{\vphi_\alp}(x)\ol{\check{\vphi_\beta}(x)}dx =  
\pair{F(\check{\vphi_\alp})}{F(\check{\vphi_\beta})}, \qquad (\alp, \beta \in \Lam_2).
\end{eqnarray}
Thus we have the Plancherel formula as following.
\begin{thm}$($Plancherel Formula$)$\label{th: size2 Plancherel}
Define the inner product on $\C[q^{x},q^{-x}, q^y+q^{-y}]$ by
\begin{eqnarray}
\pair{f}{g} = \int_T dx \int_U f\ol{g}w(y)dy,
\end{eqnarray}
where
\begin{eqnarray*}
&&
x = \dfrac{z_2+z_1}{2}, \;  y = \dfrac{z_2-z_1}{2}, \; T = \set{x = \sqrt{-1}t}{0 \leq t \leq 2\pi\log q},\\
&&
U \mbox{ is the set indicated in {\rm Figure $2$}},\;  
\int_{T\times U}dxdy = \dfrac{(1+q^{-2})^2}{2(1-q^{-1})},\\
&&
w(y) = \dfrac{1-q^{2y}}{(1-q^{2y+1})(1-q^{2y-2})}\cdot \dfrac{1-q^{-2y}}{(1-q^{-2y+1})(1-q^{-2y-2})}.
\end{eqnarray*}
Then, for any $\vphi, \psi \in \SKX$, the following identity holds:
 \begin{eqnarray}
\int_X \vphi(x)\ol{\psi(x)}dx = \pair{F(\vphi)}{F(\psi)}.
\end{eqnarray}
\end{thm}

\medskip
As a corollary of Plancherel formula, we have

\begin{cor} $($Inversion Formula$)$
For any $\vphi \in \SKX$, the following identity holds: 
\begin{eqnarray*}
\vphi(x) = \dfrac{1}{v(K\cdot x)} \pair{F(\vphi)}{F(ch_x)},  \quad (x \in X),
\end{eqnarray*}
where $ch_x$ is the characteristic function of $K\cdot x$ in $\SKX$.
\end{cor}

\vspace{2cm}
\noindent
Acknowledgement:  The author expresses her sincere gratitude to Satoshi Murai for his advice to use Macaulay2(\cite{Macaulay}) and to Yasushi Komori for the helpful discussion on orthogonal polynomials and functions. The author is grateful to the referee for careful reading and kind comments.  

\vspace{1cm}
\bibliographystyle{amsalpha}

\begin{thebibliography}{HS}
%
%
%
%

%
%
\bibitem[H1]{H1JapJM} Y.~Hironaka:
Spherical function of hermitian and symmetric forms I, {\it Japan. J. Math. }{\bf 14}(1988), 203 -- 223.

\bibitem[H2]{JMSJ} Y.~Hironaka:
Spherical function and local densities on hermitian forms,
{\it J. Math. Soc. Japan} {\bf 51}(1999), 553 -- 581.

\bibitem[H3]{JNT} Y.~Hironaka:
Local zeta functions on hermitian forms and its application to local densities, 
{\it J. Number Theory }{\bf 71}(1998), 40 -- 64.

\bibitem[H4]{French} Y.~Hironaka:
Spherical functions on $p$-adic homogeneous spaces,\\
``Algebraic and Analytic Aspects of Zeta Functions and 
$L$-functions -- Lectures at the French-Japanese Winter 
School (Miura, 2008)-- ", 
{\it MSJ Memoirs} {\bf 21}(2010), 50 -- 72. 


\bibitem[H5]{Heven} 
Y.~Hironaka:
{Harmonic analysis on the space of $p$-adic unitary hermitian matrices, mainly for dyadic case},
{\it Tokyo J. Math. }{\bf 40}(2017), 517 -- 564.

%
%
\bibitem[HK1]{HK1} Y.~Hironaka and Y.~Komori:
{Spherical functions on the space of $p$-adic unitary hermitian matrices},
{\it Int.~J.~Number Theory }{\bf 10}(2014). 513 -- 558. 
 
\bibitem[HK2]{HK2} Y.~Hironaka and Y.~Komori:
{Spherical functions on the space of $p$-adic unitary hermitian matrices II, the case of odd size},
{\it Commentarii Math. Univ. Sancti Pauli} {\bf 63}(2014), 47 -- 78. 
%

\bibitem[HS1]{HS1} Y.~Hironaka and F.~Sato:
Spherical functions and local densities of alternating
   forms,   
   {\it American Journal of Mathematics } {\bf 110}(1988), 473 -- 512.

%
%
%
%
\bibitem[HS2]{HS2} Y.~Hironaka and F.~Sato:
Local densities of alternating forms, 
{\it Journal of Number Theory } {\bf 33}(1989), 32 -- 52.

\bibitem[Jac]{Jac} R.~Jacobowitz:
Hermitian forms over local fields, {\it Amer. J. Math. }{\bf 84}(1962), 12 -- 22.
%

\bibitem[Ki1]{Kitaoka} Y.~Kitaoka:
Representations of quadratic forms and their application to Selberg's zeta functions, 
{\it Nagoya Math. J. }{\bf 63}(1976), 153 -- 162.

\bibitem[Ki2]{Kita2} Y.~Kitaoka:
A note on local densities of quadratic forms, 
{\it Nagoya Math. J. }{\bf 92}(1983), 145 -- 152.


\bibitem[KS]{KS96} R.~Koekoek and R.~F.~Swarttouw: 
The Askey-scheme of hypergeometric orthogonal polynomials and its $q$-analogue,
  arXiv math/9602214.

\bibitem[M1]{Mac} I. G. Macdonald:
{\it Symmetric Functions and Hall Polynpomials}, Oxford Science Publ., 1979.


\bibitem[M2]{Mac2} I. G. Macdonald:
Orthogonal polynomials associated with root systems, {\it S{\' e}minaire Lothanringien de Combinatoire} 
{\bf 45}(2000). Article B45a.

\bibitem[Mac2]{Macaulay}  D. Grayson and M. Stillman: Macaulay2, a software system for research in algebraic
geometry. Available at http://www.math.uiuc.edu/Macaulay2/.

\bibitem[Mar]{Martin} Kimball Martin: (Quaternion) Algebras in Number Theory, Course Notes, 2017,
 http://www2.math.ou.edu/~kmartin/quaint/ .

\bibitem[Oh] {OY} Yasuhiro Ohtaka: Spherical functions on Quaternion hermitian matrices (Master Thesis, in Japanese), 2004.1.

\bibitem[Re]{MO} I.~Reiner: {\it Maximal Orders}, Academic Press, 1975; reissued as London Mathematical Society Monograph New Series {\bf 28}, 2003.
%
%
%
%
%
%
%
%
\end{thebibliography}


\vspace{1cm}
\begin{flushleft}
Yumiko Hironaka\\
Department of Mathematics\\
Faculty of Education and Integrated Sciences \\
Waseda University\\
Nishi-Waseda, Tokyo, 169-8050, JAPAN\\
e-mail: hironaka@waseda.jp
\end{flushleft}

\end{document}